\documentclass[aop,MSNbibl,seceqn,secthm,dvips]{arximspdf}
\usepackage{graphicx}

\doi{10.1214/10-AOP625}
\volume{40}
\issue{2}
\pubyear{2012}
\firstpage{714}
\lastpage{742}

\makeatletter

\let\epsilon\varepsilon

\newcommand{\gfrac}[2]{#1/#2}
\newcommand{\eqref}[1]{(\ref{#1})}
\newcommand{\cA}{\mathcal{A}}
\newcommand{\cB}{\mathcal{B}}
\newcommand{\cC}{\mathcal{C}}
\newcommand{\cD}{\mathcal{D}}
\newcommand{\cE}{\mathcal{E}}
\newcommand{\cF}{\mathcal{F}}
\newcommand{\cH}{\mathcal{H}}
\newcommand{\cL}{\mathcal{L}}
\newcommand{\cN}{\mathcal{N}}
\newcommand{\cQ}{\mathcal{Q}}
\newcommand{\cT}{\mathcal{T}}
\newcommand{\cY}{\mathcal{Y}}
\newcommand{\frD}{\mathfrak{D}}
\newcommand{\frT}{\mathfrak{T}}
\newcommand{\frW}{\mathfrak{W}}
\newcommand{\frl}{\mathfrak{l}}
\newcommand{\frq}{\mathfrak{q}}
\newcommand{\frr}{\mathfrak{r}}
\newcommand{\frs}{\mathfrak{s}}
\newcommand{\bbB}{{\mathbb{B}}}
\newcommand{\bbC}{{\mathbb{C}}}
\newcommand{\bbE}{{\mathbb{E}}}
\newcommand{\bbN}{{\mathbb{N}}}
\newcommand{\bbP}{{\mathbb{P}}}
\newcommand{\bbR}{{\mathbb{R}}}
\newcommand{\bbS}{{\mathbb{S}}}
\newcommand{\bbZ}{{\mathbb{Z}}}
\newcommand{\gb}{\beta}
\newcommand{\gd}{\delta}
\newcommand{\go}{\omega}
\newcommand{\gl}{\lambda}
\newcommand{\sfw}{{\mathsf w}}
\newcommand{\sft}{{\mathsf t}}
\newcommand{\sfe}{{\mathsf e}}
\newcommand{\sfU}{{\mathsf U}}
\newcommand{\sfV}{{\mathsf V}}
\newcommand{\sfX}{{\mathsf X}}
\newcommand{\sfY}{{\mathsf Y}}
\newcommand{\sfZ}{{\mathsf Z}}
\newcommand{\sfP}{{\mathsf P}}
\newcommand{\bfT}{{\mathbf{T}}}
\newcommand{\bfl}{{\mathbf{l}}}
\newcommand{\bfq}{{\mathbf{q}}}
\newcommand{\bfr}{{\mathbf{r}}}
\newcommand{\1}{\mathbf{1}}
\newcommand{\smallo}{o}
\newcommand{\so}{\smallo(1)}
\newcommand{\df}{\stackrel{\Delta}{=}}
\newcommand{\leqs}{\lesssim} % Less up to a constant
\newcommand{\geqs}{\gtrsim} % Greater up to a constant
\newtheorem{lem}[thm]{Lemma} % Numbered along with thm
\newtheorem{prop}[thm]{Proposition} % Numbered along with thm
\newcommand{\cone}{\cY_\delta}
\newcommand{\Prw}{\mathrm{P}_{\mathrm{RW}}}
\newcommand{\Prwt}{\Prw^\otimes}
\newcommand{\PRWS}{\widehat{\mathrm{P}}}
\newcommand{\ERWS}{\widehat{\mathrm{E}}}
\newcommand{\Tr}{{\mathsf{Tr}}}
\newcommand{\Pd}[1]{\sfP_{#1}^d}
\newcommand{\Pslim}{\operatorname{\mathbb{P}^*\mbox{-}\lim}}
\newcommand{\tendpt}{\pi^{\perp}}
\newcommand{\lendpt}{\pi^{\parallel}}
\newcommand{\dd}{d} % a straight d for differentials
\renewcommand{\emptyset}{\varnothing}
\newcommand{\sumtwo}[2]{\mathop{\sum_{#1}}_{ #2}} % sum with 2 lines
\newcommand{\prodtwo}[2]{\mathop{\prod_{#1}}_{ #2}} % product 2
\newproclaim{rem}[thm]{Remark}

\newproclaim{ass}{Assumption}

\newtheorem{bigthm}{Theorem}

\makeatother

\begin{document}
\begin{frontmatter}

\title{Crossing random walks and stretched polymers at~weak disorder}
\runtitle{Crossing random walks and stretched polymers}

\begin{aug}
\author[A]{\fnms{Dmitry} \snm{Ioffe}\ead[label=e1]{ieioffe@ie.technion.ac.il}\thanksref{aut1}}
\thankstext{aut1}{Supported in part by the Israeli Science Foundation Grant No 817/09.}
and
\author[B]{\fnms{Yvan} \snm{Velenik}\corref{}\ead[label=e2]{Yvan.Velenik@unige.ch}\thanksref{aut2}}
\thankstext{aut2}{Supported in part by the Swiss National Science Foundation.}

\runauthor{D. Ioffe and Y. Velenik}
\affiliation{Technion and Universit\'{e} de Gen\`{e}ve}
\address[A]{Faculty of Industrial Engineering\\
Technion\\
 Haifa 32000\\
Israel\\
\printead{e1}} %adresu isvedimo komanda gale!
\address[B]{Section de Math\'{e}matiques\\
Universit\'{e} de Gen\`{e}ve\\
2--4 rue du Li\`{e}vre\\
1211 Gen\`{e}ve 4\\
Switzerland\\
\printead{e2}}
\end{aug}

% HISTORY:
\received{\smonth{3} \syear{2010}}
\revised{\smonth{9} \syear{2010}}

% ABSTRACT
%
\begin{abstract}
We consider a model of a polymer in $\mathbb{Z}^{d+1}$, constrained to
join $0$ and a hyperplane at distance $N$. The polymer is subject to a
quenched nonnegative random environment. Alternatively, the model
describes crossing random walks in a random potential (see Zerner
[\textit{Ann Appl. Probab.} \textbf{8} (1998) 246--280] or Chapter~5 of
Sznitman [\textit{Brownian Motion, Obstacles and Random Media} (1998)
Springer] for the original Brownian motion formulation). It was
recently shown [\textit{Ann. Probab.} \textbf{36} (2008) 1528--1583;
\textit{Probab. Theory Related Fields} \textbf{143} (2009) 615--642]
that, in such a setting, the quenched and annealed free energies
coincide in the limit $N\to\infty$, when $d\geq3$ and the temperature
is sufficiently high. We first strengthen this result by proving that,
under somewhat weaker assumptions on the distribution of disorder
which, in particular, enable a small probability of traps, the ratio of
quenched and annealed partition functions actually converges. We then
conclude that, in this case, the polymer obeys a diffusive scaling,
with the same diffusivity constant as the annealed model.
\end{abstract}

% KEYWORDS
%

\end{frontmatter}

\setcounter{footnote}{2}

%s1 ###
\section{Notation and results}
For simplicity\footnote{We emphasize that our techniques can deal in
the same way with any finite-range step distribution. Similarly, the
particular geometric setting used, with the arrival hyperplane
orthogonal to some lattice direction, can easily be generalized.} we
shall consider stretched polymers which are represented by
nearest-neighbor paths on $\bbZ^{d+1}$. Due to the presence of a
preferred direction, it is convenient to decompose ${\mathsf x}\in\bbZ
^{d+1}$ into transverse and longitudinal parts: ${\mathsf x}=({\mathsf
x}^{\perp},{\mathsf x}^{\parallel})$ with ${\mathsf x}^{\perp}\in\bbZ^d$ and
${\mathsf x}^{\parallel}\in\bbZ$. Given $N\in
\bbN$, we define
\[
\cH_N^- \df \{{\mathsf x}\in\bbZ^{d+1}  \dvtx  {\mathsf
x}^{\parallel} < N  \}
\]
and its outer vertex boundary $\cL_N \df\partial\cH_N^-$. We shall
consider the family~$\cD_N$ of nearest-neighbor paths from the origin
$0$ to $\cL_N$. The name \textit{stretched} stipulates that although the
second endpoint of $\gamma\in\cD_N$ is constrained to lie on $\cL
_N$, there are no other restrictions on the geometry of polymers, which
can bend and self-intersect. In the Brownian version of this
problem \cite{Sznitman}, an alternative designation often used in the
literature is \textit{crossing Brownian motion}.

The weight $W_{\gl,\beta}^\go(\gamma)$ of a polymer $\gamma= (
\gamma(0), \dots,\gamma(n))\in\cD_N$ is given by
%
%e1.1 ###
\begin{equation}
\label{eq:qweights}
W_{\gl,\beta}^\go(\gamma) \df
\exp \Biggl\{ -\lambda n - \beta\sum_{l=1}^{n} V^\go(\gamma
(l)) \Biggr\} .
\end{equation}
Here $\gl>\gl_0\df\log(2 d+2)$, $\beta>0$ and the random
environment $\{V^\go(x)\}_{x\in\bbZ^{d+1}}$, $\omega\in\Omega$,
is assumed to be i.i.d., $ V^\go(x )\stackrel{\mathrm{d}}{\sim} V$, and
such that:

\renewcommand{\theass}{(\Alph{ass})}
\begin{ass}\label{assA}
$0\in\operatorname{supp} (V )\subseteq[0,\infty]$ and $p\df\bbP(V =\infty
)$ is
sufficiently small.
\end{ass}

That the potential $V$ be bounded below is essential, since it guarantees
ballistic behavior (spatial extension) of stretched polymers.

%attractivity of the
The condition on the smallness of $p$ is also essential, since it
guarantees that we never meet situations when $\{{\mathsf x}: V^\go
({\mathsf x})<\infty\}$ does not percolate.
On the other hand, the condition $\inf\operatorname{supp} (V) = 0$ is just a
normalization.

The corresponding quenched and annealed partition functions are defined as
\[
\frD^\go_N
= \frD^\go_N (\gl,\gb)
\df\sum_{\gamma\in\cD_N}W_{\gl,\beta}^\go(\gamma) \quad \mbox{and}
\quad
{\mathbf D}_N \df\bbE\frD^\go_N .
\]
Note that the annealed potential is always attractive: For any pair of
paths~$\gamma_1$ and $\gamma_2$,
%
%e1.2 ###
\begin{equation}
\label{eq-attract}
\bbE ( W_{\gl,\beta}^\go(\gamma_1) W_{\gl,\beta}^\go
(\gamma_2)  ) \geq\bbE ( W_{\gl,\beta}^\go(\gamma_1)
 ) \bbE ( W_{\gl,\beta}^\go(\gamma_2)  ).
\end{equation}
(This can be most easily deduced from the fact that decreasing
functions on~$\bbR$ are always positively correlated.)

It has recently been proved by Flury \cite{Flury} (under the
additional \mbox{assumption} that $\bbE V^{d+1}\!<\!\infty$), and then reproved
by Zygouras \cite{Zygouras} [for arbitrary \mbox{directions},
under the additional assumption that
$\operatorname{supp} (V)$ be bounded] that, in four and higher dimensions (i.e.,
for $d\geq3$ in our notation) and for any $\gl>\gl_0$, the annealed
and quenched free energies coincide when $\beta$ is small enough.
Namely, for all $\beta$ sufficiently small, there exists $\xi=\xi
(\gl,\beta) >0$ such that
%
%e1.3 ###
\begin{equation}
\label{Flury}
-\lim_{N\to\infty}\frac{1}{N}\log\frD^\go_N = \xi=
-\lim_{N\to\infty}\frac{1}{N}\log{\mathbf D}_N.
\end{equation}
This is an important result: In sharp contrast with models of directed
polymers, the model of stretched\vadjust{\goodbreak} polymers does not have an immediate
underlying martingale structure, and this makes it necessary to find
different (and arguably more intrinsic) approaches to its analysis. The
condition $\bbE V^{d+1} <\infty$, under which \eqref{Flury} was
derived, is inherited from \cite{Zerner}, where it was shown to be sufficient
to guarantee the existence of the quenched free energy, that is, the
left-most limit
in \eqref{Flury}.

In the sequel, we shall prove the following sharp version of \eqref{Flury}:
Let $\operatorname{Cl}_\infty(V)$ be the (unique) infinite connected cluster of
sites ${\mathsf x}$ with $V ({\mathsf x})<\infty$. Under
Assumption \ref{assA}, such
a cluster $\bbP$-a.s.\ exists and is unique.

\renewcommand{\thebigthm}{\Alph{bigthm}}
\begin{bigthm}
\label{Thm:limit}
Let $d\geq3$. Then, for every $\lambda>\lambda_0$, there exists
$\beta_0 =\beta_0 (\gl,d)$ and
$p_\infty>0$, such that, if  Assumption \ref{assA} holds
with $p\leq p_\infty$,
then, for every $\beta\in[0 ,\beta_0 )$, the limit
%
%e1.4 ###
\begin{equation}
\label{eq:limit}
{\mathfrak{d}^\go} \df\lim_{N\to\infty} \frac{\frD
^\go_N}{{\mathbf D}_N}
%  \in  (0,\infty)
\end{equation}
exists $\bbP$-a.s.\ and in $L^2 (\Omega)$. In particular, the
quenched free energy\break
$-\lim_{N\to\infty}\frac{1}{N}\times\log\frD^\go_N$ is well defined,
and \eqref{Flury} holds. Furthermore, $\mathfrak{d}^\go
>0$ $\bbP$-a.s.\ on the
event $\{0\in\operatorname{Cl}_\infty(V)\}$.
\end{bigthm}

Our work was inspired by \cite{Flury,Zygouras}; however, our proof of
Theorem \ref{Thm:limit} does not rely on the results therein. In
particular, in addition to strengthening their conclusion, Theorem \ref
{Thm:limit} lifts some of the restrictions imposed on
the potential $V$ in these works.
In fact, under our assumptions, which do not impose any moment conditions
on the distribution of $V$ and even enable a~small probability of traps,
the existence of the quenched free energy needs a justification:
as we have already mentioned, the corresponding existence
results in \cite{Zerner}, which is a reference work for both
\cite{Flury} and \cite{Zygouras}, have been established under the
additional assumption
$\bbE V^{d+1} <\infty$.

Our second result confirms the prediction that stretched polymers
should be diffusive at weak disorder: On the event
$0\in\operatorname{Cl}_\infty(V)$, the random\break weights~\eqref{eq:qweights} induce a
(random) probability distribution $\mu_N^\go$ on $\cD_N$. For a
polymer $\gamma= (\gamma(0), \dots,\gamma(n))\in\cD_N$, we define
$\tendpt(\gamma)$ as the ($\bbZ^d$-valued) transverse component of
its endpoint, and $\lendpt(\gamma)=N$ as its longitudinal component,
so that $\gamma(n) = (\tendpt(\gamma), \lendpt(\gamma))$.

\begin{bigthm}
\label{Thm:diffusive}
Let $d\geq3$. Then, for every $\lambda>\lambda_0$, there exist $\hat
\beta_0 =\hat\beta_0 (\gl,d)$ and $\hat p_\infty>0$ such that, if
Assumption \ref{assA} holds with $p \leq\hat p_\infty$, then, for
every $\beta\in[0, \hat\beta_0)$, the distribution of $\tendpt$
displays diffusive scaling with a nonrandom nondegenerate diffusivity
matrix $\Sigma$ and, accordingly, a positive diffusivity
constant
$\sigma^2 = \sigma^2 (\beta,\gl) \df\Tr(\Sigma) > 0$. Namely,
define $\bbP^* (\cdot) \df\bbP(\cdot| 0\in\operatorname{Cl}_\infty(V) )$. Then,
%
%e1.5 ###
\begin{equation}
\label{eq:diffusive}
\Pslim\limits_{N\to\infty}\mu_N^\go \biggl(\frac{|\tendpt(\gamma)|^2}N
 \biggr)
= \sigma^2 ,\vadjust{\goodbreak}
\end{equation}
where $\Pslim$ denotes convergence in $\bbP^*$-probability.
Furthermore, for any bounded continuous function $f$ on $\bbR^d$,
\begin{eqnarray}
\label{eq:flim}
&&\operatorname{\Pslim}\limits_{N\to\infty} \sum_{{\mathsf x}\in\bbZ^d}\mu_N^\go\bigl(
\tendpt
(\gamma) = {\mathsf x}\bigr)f\biggl(\frac{{\mathsf
x}}{\sqrt{N}}\biggr)\nonumber\\[-9pt]\\[-9pt]
&&\qquad = \frac{1}{\sqrt{\operatorname{det}(2\pi\Sigma)}}\int_{\bbR^d} f({\mathsf x})
 e^{-\gfrac12 (\Sigma^{-1}{\mathsf x},{\mathsf x})}
\,\dd{\mathsf x}.\nonumber
\end{eqnarray}
$\Sigma$ and $\sigma^2$ above are precisely the diffusion matrix and
the diffusivity constant of the corresponding annealed polymer model;
see \eqref{eq:AnnDiff} below.\vspace*{-1pt}
\end{bigthm}

We expect both \eqref{eq:diffusive} and \eqref{eq:flim} to hold not
only in $\bbP^*$-probability, but also in $L^2(\Omega)$ and $\bbP^*$-a.s.\vspace*{-1pt}

%s1.1 ###
\subsection{Some open problems}
In this subsection, we briefly discuss some points that are left
untouched in the present work.\vspace*{-1pt}

\subsubsection*{Stronger modes of convergence}
As already mentioned above, we expect our diffusivity results to hold
also a.s.\ in the environment and in $L^2(\Omega)$. Such results are
known in the directed case, as a consequence of the much simpler
martingale structure \cite{Bolthausen1989}. Furthermore, we expect the
$\bbP^*$-a.s.\ validity of a~local CLT, or equivalently, of a (random)
Ornstein--Zernike-type formula for long-range quenched connections; see
the discussion at the end of Section~\ref{sec:ProofThmDiffusive}.\vspace*{-1pt}

\subsubsection*{Invariance principle}
Once equipped with a local CLT and thanks to our good control on the
path geometry, it should be mostly straightforward to obtain a full
invariance principle for the path.\vspace*{-1pt}

\subsubsection*{``Real'' stretched polymer}
In the present work, we have focused on ensem\-bles of paths of
``point-to-plane'' type (the set $\cD_N$). It would be physi\-cally
quite interesting to analyze also the case of fixed-length polymers,
stretched by an external force (notice that in the directed case there
is no difference between
``point-to-plane'' and ``fixed-length'' scenarios); in particular, it
would be interesting to obtain a local limit theorem for the free
endpoint. Such questions have been investigated in the annealed setting
in our previous work \cite{IV-annealed}. In the quenched setting
coincidence of Lyapunov exponents (under the
additional $\bbE V^{d+1} <\infty$ assumption) has been established in
\cite{Flury}.\vspace*{-1pt}

\subsubsection*{Nonperturbative proof}
Our results are only valid at very high temperatures. It would be quite
interesting (and probably challenging) to push them to the full
weak-disorder regime. Results of that type have been obtained in the
directed case~\cite{CometsYoshida2006}.\vspace*{-1pt}

\subsubsection*{Strong disorder}
We only consider the weak disorder case here. Obtaining some
information on the behavior\vadjust{\goodbreak} of typical paths in the strong disorder
regime would also be quite interesting, and is the subject of some work
in progress.
See \cite{CarmonaHu2006} for such results, in the full strong disorder
regime, in the directed case.

%s1.2 ###
\subsection{A remark on notational conventions}
Given two sequences $\{ a_n ( w)\}$ and $\{ b_n (w )\}$ of positive
real numbers indexed by $w$ from some set of parameters $\frW_n$, we
say that $a_n ( w)\leqs b_n (w )$, if
\[
\limsup_{n\to\infty}\frac{a_n (w)}{b_n (w)}  <  \infty,
\]
uniformly in $w\in\frW_n$.

Given ${\mathsf z},\sfw\in\bbC^{d+1}$, we use
\[
({\mathsf z},\sfw)_{d+1} \df\sum_{i=1}^{d+1}{\mathsf z}_i\bar
\sfw_i \quad \mbox{and}\quad
({\mathsf z},\sfw)_{d} \df\sum_{i=1}^{d}{\mathsf z}_i\bar\sfw_i.
\]
With a slight abuse of notation, we shall also write $(z ,\sfw)
_{d}$ for the same expression with $z\in\bbC^d$.

%s2 ###
\section{Convergence of partition functions}
%s2.1 ###
\subsection{\texorpdfstring{Irreducible decomposition of paths $\gamma\in\cD_N$}
{Irreducible decomposition of paths gamma in D_N}}

Given $\delta>0$, we define a positive cone along the $\vec{\sfe}\df
\vec{\sfe}_{d+1}$-direction by
\[
\cY_\delta\df \{ {\mathsf x}\in\bbR^{d+1}  \dvtx  \|{\mathsf
x}^{\perp}\| <\delta x^{
\parallel
} \} ,
\]
where $\|\cdot\|$ denotes the Euclidean norm.
We say that a trajectory $\gamma= (\gamma(0) ,\ldots, \gamma(n))$
of length $|\gamma| = n$ is \emph{cone-confined} if
\[
\gamma\subseteq\bigl(\gamma(0)+\cY_\delta\bigr)\cap\bigl(\gamma(n)-
\cY_\delta\bigr).
\]
Although paths $\gamma\in\cD_N$ always satisfy $0 =(\gamma(0), \vec
{\sfe})_{d+1} < (\gamma(n), \vec{\sfe})_{d+1}$, evi\-dently not all
of them are cone-confined.
For $1\leq k <n =|\gamma|$, let us say that~$\gamma(k)$ is a \emph
{cone-point} of $\gamma$ if
\[
(\gamma(0), \vec{\sfe})_{d+1} < (\gamma(k), \vec{\sfe})_{d+1} <
(\gamma(n), \vec{\sfe})_{d+1} ,
\]
and, in addition, if
\[
\gamma\subseteq\bigl(\gamma(k)-\cY_\delta\bigr)\cup\bigl(\gamma(k)+
\cY_\delta\bigr).
\]
We say that a trajectory $\gamma$ is \emph{irreducible} if it
contains less than two cone-points.
We say that it is \emph{strongly irreducible} if it does not contain
cone-points at all.

%f1 ###
\begin{figure}

\includegraphics{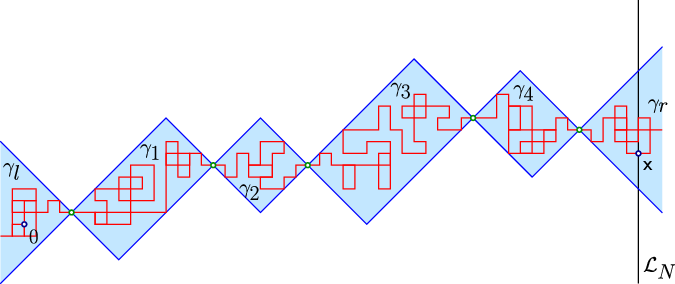}

\caption{The decomposition of a path $\gamma\in\cD_N$ into a
concatenation of strongly irreducible pieces.}
\label{fig:decomp}
\end{figure}

The following mass-separation property of irreducible trajectories,
proved in \cite{IV-annealed}, is crucial to our analysis: There exists
$\nu>0$ such that, for all $N$ large enough,
%
%e2.1 ###
\begin{equation}\label{eq-SepMass}
\frac1{\mathbf{D}_N}\mathop{\sum_{\gamma\in\cD_N}}_{\mathrm{irreducible}}
\bbE W_{\gl,\beta}^\go(\gamma) \leq e^{-\nu N}.\vadjust{\goodbreak}
\end{equation}
On the other hand, reducible trajectories are unambiguously represented
as concatenation of strongly irreducible pieces (as induced by the
collection of all the cone-points of $\gamma$; see Figure \ref{fig:decomp}),
%
%e2.2 ###
\begin{equation}\label{eq:gamma-split}
\gamma= \gamma_l\cup\gamma_1\cup\cdots\cup\gamma_n\cup\gamma_r .
\end{equation}
By construction, $\gamma_1 ,\dots,\gamma_n$ above are also
cone-confined, and so is their concatenation $\gamma_1\cup\cdots\cup
\gamma_n$. Thus, \eqref{eq-SepMass} and \eqref{eq:gamma-split}
suggest that the asymptotics of $\mathbf{D}_N$ and $\frD_N^\go$
should be closely related to the asymptotics of the corresponding
cone-confined quantities. This intuition turns out to be correct.

Let $\cT_N$ be the family of all cone-confined trajectories from $0$
to $\cL_N$. Set
\[
\frT_N^\go(\gl,\beta) \df
\sum_{\gamma\in\cT_N}W_{\gl,\beta}^\go(\gamma)
\quad \mbox{and}\quad
{\mathbf T}_N \df\bbE\frT^\go_N .
\]
The following statement as well as the understanding one needs to
develop for its proof are crucial: In the notation and under the
conditions of Theorem~\ref{Thm:limit}, for every $\beta\in[0 ,\beta
_0 )$, the limit
%
%e2.3 ###
\begin{equation}
\label{eq:limitT}
\lim_{N\to\infty} \frac{\frT^\go_N}{{\mathbf T}_N}
% \in(0,\infty)
\end{equation}
exists $\bbP$-a.s.\ and in $L^2 (\Omega)$. For a while we shall focus
on the ensembles of cone-confined trajectories and on proving \eqref
{eq:limitT}. We shall return to $\cD_N$ and prove the full
statement \eqref{eq:limit} only in Section \ref{sub:DN}.

\subsubsection*{Notation for scaled quantities}
Recall the definition of the Lyapunov exponent $\xi$ in \eqref
{Flury}. Given $N\geq1$ and $\gamma\in\cT_N$, we define the scaled
random path weights
\[
w_{\gl,\beta}^\go(\gamma) \df e^{N\xi} W_{\gl,\beta}^\go
(\gamma) .\vadjust{\goodbreak}
\]
For ${\mathsf x}\in\cL_N$, we define
\begin{eqnarray} \label{eq:txomega}
\quad \mathfrak{t}_{\mathsf x}^\go&\df&
\sumtwo{\gamma:0\to{\mathsf x}}{\gamma\in\cT_N}
w_{\gl,\beta}^\go(\gamma),\qquad
\frq_{\mathsf x}^\go\df\sumtwo{\gamma:0\to{\mathsf x}}{\gamma\in
\cT_N^0}
w_{\gl,\beta}^\go(\gamma)
\quad \mbox{and}\nonumber\\[-8pt]\\[-8pt]
{\mathbf{ t}} _{\mathsf x}&\df&
\bbE\mathfrak{t}_{\mathsf x}^\go,\qquad
\bfq_{\mathsf x}\df
\bbE\frq_{\mathsf x}^\go,\nonumber
\end{eqnarray}
where $\cT_N^0$ denotes the set of all strongly irreducible $\gamma
\in\cT_N$. Similarly, we define
\[%\label{eq:qxomega}
\mathfrak{t}_N^\go\df\sum_{{\mathsf x}\in\cL
_N}\mathfrak{t}_{\mathsf x}^\go,\qquad
\frq_N^\go\df\sum_{{\mathsf x}\in\cL_N}\frq_{\mathsf x}^\go
\quad \mbox{and}\quad
{\mathbf{ t}} _N \df\bbE\mathfrak{t}_N^\go,\qquad
\bfq_N \df\bbE\frq_N^\go.
\]
We also set $\mathfrak{t}^\go_0={\mathbf{ t}}
_0\df1$.

%s2.2 ###
\subsection{Renewal analysis of annealed partition function ${\mathbf T}_N$}

With the above notation, the sequence $\{{\mathbf{ t}}
_N\}$ satisfies the renewal relation
%
%e2.4 ###
\begin{equation}\label{eq:annrenewal}
{\mathbf{ t}} _0 = 1\quad \mbox{and}\quad
{\mathbf{ t}} _N = \sum_{M=0}^{N-1} {\mathbf{ t}} _M   \bfq_{N-M},\qquad
  N\geq1 .
\end{equation}
We fix $\lambda$ and $\beta$ and set
%
%e2.5 ###
\begin{equation}
\label{eq:mu}
\mu= \mu(\lambda,\beta)
\df\sum_{M\geq1} M \bfq_M .
\end{equation}
Note that the above series converges since, by our basic
mass-separation estimate for annealed quantities \eqref{eq-SepMass},
%
%e2.6 ###
\begin{equation}\label{eq:massgap}
\bfq_M \leq e^{-\nu M}  e^{M\xi}\mathbf{D}_M \leq e^{-\nu M},
\end{equation}
where we used the fact that $\mathbf{D}_M \leq e^{- M\xi}$,
which follows from subadditivity.
\begin{lem}
For any $\beta\geq0$ and $\lambda>\lambda_0$,
%
%e2.7 ###
\begin{equation}\label{eq:convtN}
\lim_{N\to\infty}  e^{N\xi} \bfT_N
= \lim_{N\to\infty} {\mathbf{ t}} _N
= \frac1{\mu(\lambda,\beta)}.
\end{equation}
Moreover, the convergence in \eqref{eq:convtN} is exponentially fast.
\end{lem}

\begin{pf}
This is a standard renewal argument which we shall briefly sketch for
completeness. As a consequence of our scaling and the mass separation
property \eqref{eq-SepMass}, the radius of convergence of the
generating function
\[
\hat{\mathbf{ t}} (u)
\df\sum_{N\geq0} u^N{\mathbf{ t}} _N
\]
is equal to $1$ (see Section 3.3.6 in \cite{IV-annealed} for
details). On the other hand, it follows from \eqref{eq:massgap} that
the irreducible generating function
\[
\hat\bfq(u) \df\sum_{N \geq1} u^N\bfq_N
\]
has radius of convergence at least $1+\nu$. This implies, via standard
arguments based on \eqref{eq:annrenewal}, that $\hat\bfq(1)=1$. Of
course, $\mu= \hat\bfq^\prime(1)$. Fix $\rho\in(0,1)$. By
Cauchy's formula,
\begin{eqnarray*}
{\mathbf{ t}} _N - \frac1{\mu}
&=& \frac1{2\pi i}\int_{\bbS_\rho}
 \biggl\{
\frac{\dd u}{u^{N+1}(1-\hat\bfq(u))} - \frac{\dd u}{u^{N+1}(1-u)\mu}
 \biggr\} \\
& =&
\frac1{2\pi i}\int_{\bbS_\rho}\frac{\Delta(u)}{u^{N+1}}\,\dd u ,
\end{eqnarray*}
where $\bbS_\rho= \partial\bbB_\rho$ and
%
%e2.8 ###
\begin{equation}\label{eq:Delta}
\Delta(u) =
\frac{(\hat\bfq(u) - \hat\bfq(1)) - \hat\bfq^\prime(1)(u-1)}
{(\hat\bfq(u) - \hat\bfq(1)) \hat\bfq^\prime(1)(u-1)}.
\end{equation}
Since $\Delta$
is analytic on $\bbB_{1+\nu^\prime}$ for some $\nu^\prime\in
(0,\nu)$, the result follows.
\end{pf}
%

%s2.3 ###
\subsection{Complex tilts and annealed diffusivity}
\label{sub:tilts}
For $\delta$ small enough, let $\Pd{2\delta}\subset\bbC^d$ be the
complex polydisc with all $d$ radii equal to $2\delta$. By the
implicit function theorem (see, e.g., \cite{Kaup}) and in view of the
mass-gap estimate \eqref{eq:massgap}, the relations
\[%\label{eq:complextilts}
\varphi[0] = 0  \quad \mbox{and}\quad
\sum_{M\geq1}\sum_{{\mathsf x}\in\cL_M} \bfq_{\mathsf x}   e^{-M\varphi[z ] + (z ,{\mathsf x})_d}
\df\sum_{M\geq1} \bfq_M [z ] =1
\]
define a holomorphic function $\varphi\dvtx\Pd{2\delta}\to\bbC$. We
shall assume that
$\delta$ is so small that
%
%e2.9 ###
\begin{equation}\label{qMz}
|\bfq_M [z ]|\leqs e^{-\nu M/2} ,
\end{equation}
uniformly in $M\geq1$ and $z\in\Pd{2\delta}$.

The analysis of the previous subsection can be readily extended to
obtain the asymptotic
expansion of the moment generating functions
\[
{\mathbf{ t}} _0 [z] \df1  \quad \mbox{and\quad  for $N\geq
1$,}\qquad  {\mathbf{ t}} _N [z] \df\sum_{{\mathsf x}\in\cL
_N} {\mathbf{ t}} _{\mathsf x}   e^{-N\varphi[z ] + (
z,{\mathsf x})_d}
\]
for $z\in\Pd{2\delta}$. Indeed, ${\mathbf{ t}} _N [z]$
satisfies the renewal relation
\[
{\mathbf{ t}} _N [z]
= \sum_{M=0}^{N-1} {\mathbf{ t}} _M [z]  \bfq_{N-M}[z] .
\]
Furthermore, if $\delta$ is sufficiently small, then not only
does \eqref{qMz} hold, but there also
exists $\nu^\prime> 0$ such that, for all $z\in\Pd{2\delta}$,\vadjust{\goodbreak}
$u=1$ is the unique solution of the equation
\[
\sum_{M\geq1} u^M \bfq_M [z ] \df\hat\bfq[z](u) = 1 ,
\]
on $\bbB_{1+\nu^\prime}\subset\bbC$. We define $\mu[z]$ exactly
as in \eqref{eq:mu} by
\[
\mu[z] \df\sum_{N\geq1} N \bfq_N [z].
\]
Relying on \eqref{qMz}, we can choose $\gd$ so small that $\mu[\cdot
]$ is analytic and nonzero on $\Pd{2\delta}$. It then follows that
\[
\lim_{N\to\infty} {\mathbf{ t}} _N [z] = \frac{1}{\mu[z]},
\]
uniformly exponentially fast on $\Pd{2\delta}$.

The annealed diffusion matrix $\Sigma$ and the corresponding
diffusivity constant $\sigma^2$ in \eqref{eq:diffusive} are defined by
%
%e2.10 ###
\begin{equation}
\label{eq:AnnDiff}
\Sigma\df{\mathrm{D}}^2_d\varphi[0]  \quad \mbox{and}\quad
\sigma^2 \df \operatorname{Tr} (\Sigma),
\end{equation}
where $D^2_d\varphi$ denotes the Hessian of $\varphi$.
Now, since we have chosen $\delta$ sufficiently small to ensure that
$\mu[\cdot]$ is analytic and does not vanish on $\Pd{2\delta}$, the
functions $\log{\mathbf{ t}} _N [\cdot] + \log\mu[\cdot
] $ are analytic and exponentially small (in $N$) on $\Pd{2\delta}$.
In particular,
\[
\operatorname{Tr} \bigl( \mathrm{D}^2_d  ( \log{\mathbf{ t}} _N
[z] {+}
\log\mu[z]
 ) \bigr)
\]
is also exponentially small. This shows that the leading contribution
(in $N$) to the log-moment generating function $\log (
{\mathbf{ t}} _N[z]  e^{N\varphi[z]}  )$ of
$\tendpt(\gamma)$ under the induced measure is given by $N\varphi
[z]$. We have thus proved that
\begin{lem}
\label{lem:sigmaN}
\[%\label{eq:sigmaN}
 \biggl| \frac{1}{N{\mathbf{ t}} _N}\sum_{{\mathsf x}\in\cL
_N} \|{\mathsf x}^{\perp}\|^2  {
{\mathbf t}} _{\mathsf x}- \sigma^2
 \biggr|\leqs\frac1{N}.
\]
Furthermore, $\tendpt(\gamma)/\sqrt{N} \Rightarrow\cN( 0, \Sigma
)$ under the sequence of annealed polymer measures $\mu_N$.
\end{lem}

%s2.4 ###
\subsection{Multi-dimensional renewal relation for quenched partition
functions}
We continue to employ {the} notation introduced in \eqref{eq:txomega}.
It is immediate to check that the following analogs of \eqref
{eq:annrenewal} hold:
%
%e2.11 ###
\begin{equation}
\label{eq:MultiRenewal}
\mathfrak{t}_{\mathsf x}^\go= \sum_{{\mathsf y}}
{\mathfrak t}_{\mathsf y}^\go\frq^{\theta_{\mathsf y}\go}_{{\mathsf
x}-{\mathsf y}}
\quad \mbox{and}\quad
\mathfrak{t}_N^\go
= \sum_{M=0}^{N-1}  \sum_{{\mathsf x}\in\cL_M}\mathfrak
{t}_{\mathsf x}^\go\frq^{\theta_{\mathsf x}\go}_{N-M}
\end{equation}
for all ${\mathsf x}\in\cH_0^+\df\{{\mathsf x}\in\bbZ^{d+1}
\dvtx  {\mathsf x}^{\parallel} >0\}$ and $N\geq1$. Set
$\mathfrak{t}_0^\go\df1$, and define the generating functions
\[
\hat\mathfrak{t}^\go(u)  \df  \sum_{N=0}^{\infty}
u^N \mathfrak{t}_N^\go
\]
and
\[
\hat\frq^\go(u)  \df  \sum_{N=1}^{\infty} u^N \frq_N^\go.
\]
Since
$|\hat\mathfrak{t}^\go(u)|\leq
\hat\mathfrak{t}^\go(|u |)$ and $\bbE\hat
{\mathfrak t}^\go(\rho)
=
\hat{\mathbf{ t}} (\rho)$, the random generating~func-\break tion~%
$\hat\mathfrak{t}^\go(u)$~is $\bbP$-a.s.\ defined and
analytic in the interior of the unit disc $\bbB_1\subset\bbC$.
Similarly, the random generating function $\hat\frq^\go(u)$
is $\bbP$-a.s.\ analytic on $\bbB_{1+\nu}$ for some $\nu>0$.

We can rewrite  \eqref{eq:MultiRenewal} in terms of the generating
function as
\begin{eqnarray}\label{eq:MultiGF}
\hat\mathfrak{t}^\go(u)
&=& 1 + \sum_{M=0}^\infty u^M \sum_{{\mathsf x}\in\cL_M}
{\mathfrak t}_{\mathsf x}^\go\hat\frq^{\theta_{\mathsf x}\go} (u)\nonumber \\
&=& 1 + \hat\bfq(u) \sum_{M=0}^\infty u^M\sum_{{\mathsf x}\in\cL_M}
\mathfrak{t}_{\mathsf x}^\go\nonumber\\[-8pt]\\[-8pt]
&&{} + \sum_{M=0}^\infty u^M\sum_{{\mathsf x}\in\cL_M}
{\mathfrak t}_{\mathsf x}^\go
 \bigl( \hat\frq^{\theta_{\mathsf x}\go} (u) - \hat\bfq(u)
\bigr)\nonumber \\
&\df&  1 + \hat\bfq(u)\hat\mathfrak{t}^\go(u) + \hat
\Psi^\go(u) .\nonumber
\end{eqnarray}
Since $|\hat\bfq(u)| <1$ whenever $|u | =\rho<1$, we can record the
last computation as
\[
\hat\mathfrak{t}^\go(u) = \frac{1 +\hat\Psi^\go
(u)}{1-\hat\bfq(u)} .
\]
Therefore,
%
%e2.12 ###
\begin{equation}
\label{eq:tNRepresentation}
\mathfrak{t}_N^\go
= \frac1{2\pi i}\int_{\bbS_\rho}\frac{1 +\hat\Psi^\go
(u)}{(1-\hat\bfq(u))u^{N+1}}
\,\dd u ,
\end{equation}
$\bbP$-a.s.\  for all $\rho\in(0,1)$.

%s2.5 ###
\subsection{Recursion under $L^2$-weak disorder}
\label{sub:L2weak}
Equation \eqref{eq:tNRepresentation} is the starting point for proving
Theorem \ref{Thm:limit}.
In fact, we are going to develop a recursion for the limit in \eqref
{eq:limitT} whenever the conditions of the latter theorem are satisfied.

Let us decompose
\[
\mathfrak{t}_N^\go= \frac1\mu\frs_N^\go+\biggl(
\mathfrak{t}_N^\go- \frac1\mu\frs_N^\go\biggr),\vadjust{\goodbreak}
\]
where\footnote{Note that the definition does not depend on the particular
choice of $\rho\in(0,1)$.}
\[
\frs_N^\go  \df
\frac1{2\pi i}\int_{\bbS_\rho}
\frac{1 + \hat\Psi^\go(u)}{u^{N+1}(1-u)}\,\dd u
\]
and, accordingly,
%
%e2.13 ###
\begin{equation}
\label{eq:stdifference}
\mathfrak{t}_N^\go- \frac1\mu\frs_N^\go
=
\frac1{2\pi i}\int_{\bbS_\rho} \frac{(1 + \hat\Psi^\go
(u))\Delta(u)}{u^{N+1}}\,\dd u ,
\end{equation}
with $\Delta(u)$ defined in \eqref{eq:Delta}.

After examining the definition of $\hat\Psi^\go$ in \eqref
{eq:MultiGF}, we arrive at the following expression for $\frs_N^\go$:
\begin{eqnarray}
\label{eq:frsN}
\frs_N^\go&=&  \biggl[ \frac{1+\hat\Psi^\go(u)}{1-u}  \biggr]_N
= 1 + \sum_{M=0}^{N-1} \sum_{{\mathsf x}\in\cL_M}\mathfrak
{t}_{\mathsf x}^\go ( \frq^{\theta_{\mathsf x}\go}_{1,N-M}
- \bfq_{1, N-M} )
\nonumber\\[-8pt]\\[-8pt]
&=&{1} +     \sum_{{\mathsf x}\in\cH_N^-, {\mathsf y}\in\cH
_{N+1}^-}
\mathfrak{t}_{\mathsf x}^\go ( \frq^{\theta_{\mathsf x}\go
}_{{\mathsf y}-{\mathsf x}}
-\bfq_{{\mathsf y}-{\mathsf x}} ),\nonumber
\end{eqnarray}
where
\[
\frq^\go_{1,l}\df\sum_{k=1}^l\frq^\go_k \quad \mbox{and}\quad
\bfq_{1,l}\df\bbE\frq^{\go}_{1,l}
\]
and we used the standard notation $ [ \sum_{k\geq0} a_k u^k
 ]_N = a_N$ for expansion coefficients.

The following theorem is proved in Sections \ref{sub:ProofthmL21} and \ref{sub:ProofthmL22}.
\begin{thm}
\label{thm:L2}
For every $\lambda>\lambda_0$, there exist $\beta_0 =\beta_0 (\gl
,d)$ and
$p_\infty>0$, such that if Assumption \ref{assA} holds
with $p\leq p_\infty$, then, for every $\beta\in[0 ,\beta_0 )$:
\begin{longlist}[(2)]
\item[(1)] The sequence $\mathfrak{t}_N^\go- \frs_N^\go/\mu
$ converges to zero
$\bbP$-a.s.\ and in $L^2 (\Omega)$.
\item[(2)] The sequence $\frs_N^\go$
converges $\bbP$-a.s.\ and in $L^2 (\Omega)$ to
%
%e2.14 ###
\begin{equation}
\label{eq:sNslimit}
\frs^\go\df1+
\sum_{{\mathsf x}\in\cH_0^+} \mathfrak{t}_{\mathsf x}^\go
 ( \frq^{\theta_{\mathsf x}\go}_{1,\infty} - 1 ) ,
\end{equation}
the latter sum also converging in $L^2 (\Omega)$.
\end{longlist}
\end{thm}

Theorem \ref{thm:L2} implies that the limit in \eqref{eq:limitT}
indeed exists and, furthermore, that it is equal to
the random variable $\frs^\go$
\[
\lim_{N\to\infty} \frac{\frT^\go_N}{{\mathbf T}_N}
=
\lim_{N\to\infty} \frac{\mathfrak{t}^\go_N}{{\mathbf
t}_N } = \lim_{N\to\infty} \frs_N^\go= \frs^\go.
\]
Note that if $0\notin \operatorname{Cl}_\infty(V)$, then $
{\mathfrak t}_N^\go= 0$ for
all $N$ sufficiently large, say $N \geq N_0( \go)$.\vadjust{\goodbreak} Consequently, in
this case $\frs^\go$ is a difference
of two convergent series,\looseness=1
\[
\frs^\go= 1 +\sum_{{\mathsf x},{\mathsf y}\in\cH_{N_0}^-}
{\mathfrak t}^\go_{\mathsf x}
\frq^{\theta_{\mathsf x}\go}_{{\mathsf y}- {\mathsf x}} -
\sum_{{\mathsf x}\in\cH_{N_0}^-} \mathfrak{t}^\go_{\mathsf x}
= 0.
\]\looseness=0
Positivity of $\frs^\go$ [or rather of the full limit $
{\mathfrak d}^\go$ in
\eqref{eq:limit}] on the event $\{0 \in\operatorname{Cl}_\infty(V)\}$
is established in the
concluding Section \ref{sub:positivity} of the paper.

%s2.6 ###
\subsection{Relation with Sinai's representation}
Our representation \eqref{eq:frsN} can be seen as an effective random
walk version of the high-temperature expansion employed by Sinai
in \cite{Sinai}. Indeed, let ${\mathsf x}\in\cL_N$. Then
\begin{eqnarray*}
\label{eq:Sinai1}
\mathfrak{t}^\omega_{\mathsf x}
&=&
\sum_{n\geq0} \sum_{{\mathsf x}_1,\ldots,{\mathsf x}_n} \prod
_{k=0}^n \frq
^\omega_{{\mathsf x}_k,{\mathsf x}_{k+1}}\\
\nonumber
&=&
\sum_{n\geq0} \sum_{{\mathsf x}_1,\ldots,{\mathsf x}_n} \prod
_{k=0}^n \bfq
_{{\mathsf x}_k,{\mathsf x}_{k+1}} \Phi^\omega({\mathsf x}_0,\ldots
,{\mathsf x}_{n+1}),
\end{eqnarray*}
where we have set ${\mathsf x}_0=0$, ${\mathsf x}_{n+1}={\mathsf x}$, and
\[
\Phi^\omega({\mathsf x}_0,\ldots,{\mathsf x}_{n+1}) \df\prod
_{k=0}^n \frac
{\frq^\omega_{{\mathsf x}_k,{\mathsf x}_{k+1}}}{\bfq_{{\mathsf
x}_k,{\mathsf x}_{k+1}}} \df\prod_{k=0}^n  \bigl(1 + \phi^\omega
({\mathsf x}_k,{\mathsf x}_{k+1})  \bigr).
\]
Using the expansion
\[
\Phi^\omega({\mathsf x}_0,\ldots,{\mathsf x}_{n+1}) = \sum
_{A\subset\{
0,\ldots, n\}} \prod_{k\in A} \phi^\omega({\mathsf x}_k,{\mathsf x}_{k+1}),
\]
we obtain the representation
\[%\label{eq:basicRepr}
\mathfrak{t}^\omega_{\mathsf x}
= \sum_{n\geq0} \sum_{{\mathsf x}_1,\ldots,{\mathsf x}_n} \prod_{k=0}^n
\bfq_{ {\mathsf x}_{k} ,{\mathsf x}_{k+1}}
\sum_{A\subset\{0,\ldots, n\}}
\prod_{\ell\in A} \phi^\omega({\mathsf x}_\ell,{\mathsf x}_{\ell
+1}) .
\]
Given $n$, ${\mathsf x}_1 , \dots, {\mathsf x}_n$ and $\emptyset\neq
A\subset
\{0, \dots,n\}$, let us say that $({\mathsf x}_{k^*}, {\mathsf
x}_{k^*+1})$ is the last perturbed segment if $k^* = \max\{k\dvtx k\in A\}
$. Keeping the last perturbed segment fixed and resumming all the rest,
we arrive at
%
%e2.15 ###
\begin{equation}
\label{eq:Sinai}
\mathfrak{t}^\omega_{\mathsf x}= {\mathbf{ t}}
_{\mathsf x} +\sum_{{\mathsf y}
, {\mathsf z}
}\mathfrak{t}^\go_{\mathsf y}(\frq^{\theta_{\mathsf y}\go
}_{{\mathsf z}-{\mathsf y}}
-\bfq_{{\mathsf z}-{\mathsf y}}){\mathbf{ t}} _{{\mathsf
x}-{\mathsf z}} .
\end{equation}
Similarly, keeping the first perturbed segment fixed and resumming all
the rest, we arrive at
%
%e2.16 ###
\begin{equation}
\label{eq:Sinaib}
\mathfrak{t}^\omega_{\mathsf x}= {\mathbf{ t}}
_{\mathsf x} +\sum_{{\mathsf y}
, {\mathsf z}
}
{\mathbf{ t}} _{\mathsf y}(\frq^{\theta_{\mathsf y}\go
}_{{\mathsf z}-{\mathsf y}}
-\bfq_{{\mathsf z}-{\mathsf y}})\mathfrak{t}_{{\mathsf x}-{\mathsf
z}}^{\theta_{\mathsf z}\go} .
\end{equation}
It would have been possible to work directly with the above
representations of \mbox{$\mathfrak{t}$-quantities}. In fact,
Theorem \ref{thm:L2}(1) can be considered as the first step along
these lines: it enables us to substitute and control the $
{\mathfrak t}$-quantities by the more tractable $\frs$-quantities, as
appears in \eqref{eq:frsN}.\vadjust{\goodbreak}

Notice though that it is not clear how to prove the almost-sure
convergence in Theorem \ref{Thm:limit} without having recourse to
martingale arguments as developed in Section \ref{ssec-key}.\vspace*{-2pt}

%s2.7 ###
\subsection{Extension to the full $\cD_N$-ensemble}\label{sub:DN}

Let us go back to {Theorem} \ref{Thm:limit}. In view of \eqref
{eq-SepMass}, there is no loss in redefining $\cD_N$ as the set of all
\emph{reducible} paths from $0$ to $\cL_N$. Thus, any $\gamma\in\cD
_N$ automatically satisfies \eqref{eq:gamma-split}. By construction
(decomposition with respect to all cone-points), none of the
paths~$\gamma_l,\allowbreak \gamma_1, \dots, \gamma_n ,\gamma_r $ in \eqref
{eq:gamma-split} has cone-points. Recall that we use {the} notation
$\cT^0$ for cone-confined paths without cone-points. Thus, $\gamma_1,
\dots,\gamma_n \in\cT^0$.

Paths $\gamma_l = (\gamma_l (0), \dots,\gamma_l (m))$ satisfy
$\gamma_l\subseteq\gamma_l (m) -\cY_\delta$, and, similarly, paths
$\gamma_r = (\gamma_r (0), \dots,\gamma_r (k))$ satisfy $\gamma
_r\subseteq\gamma_r (0) +\cY_\delta$. We denote by $\cT^0_l$ and
$\cT^0_r$ the sets of such paths; in this way, $\cT^0 =\cT_l^0\cap
\cT_r^0$.

Following \eqref{eq:txomega}, define
\[
\frl^\go_{{\mathsf x}} \df\sumtwo{\gamma:0\mapsto{\mathsf
x}}{\gamma\in
\cT_l^0} w_{\lambda,\beta}^\go(\gamma)
\quad \mbox{and}\quad
\frr^\go_{{\mathsf x}} \df\sumtwo{\gamma:0\mapsto{\mathsf
x}}{\gamma\in
\cT_r^0} w_{\lambda,\beta}^\go(\gamma) .
\]
As usual, we denote the corresponding annealed quantities by $\bfl
_{{\mathsf x}}$ and $\bfr_{{\mathsf x}}$.
The scaled full $\cD_N$ partition function satisfies
%
%e2.17 ###
\begin{eqnarray}
\label{eq:dNDecomp}
\mathfrak{d}_N^\go&\df&
 e^{N\xi}\frD_N^\go= \sum_{{\mathsf x}\in\cL_N}\sumtwo
{\gamma
:0\mapsto{\mathsf x}}{\gamma\in\cD_N}
w_{\lambda,\beta}^\go(\gamma)
=\sum_{0 \leq M_l < M_r\leq N}  \sumtwo{{\mathsf x}\in\cL
_{M_l}}{{\mathsf y}\in\cL_{M_r}}
\frl^\go_{{\mathsf x}}\mathfrak{t}_{{\mathsf y}-{\mathsf x}}^{\theta
_{\mathsf x}\go}\frr_{N-M_{r} }^{\theta_{\mathsf y}\go} \nonumber \\[-1pt]
&=&
\sum_{0 \leq M_l < M_r\leq N}\  \sum_{{\mathsf x}\in\cL_{M_l}}
\frl^\go_{{\mathsf x}}\mathfrak{t}_{M_r -M_l}^{\theta_{\mathsf
x}\go}\bfr_{N-M_{r}}\\[-1pt]
&&{} +
\sum_{0 \leq M_l < M_r\leq N}\  \sumtwo{{\mathsf x}\in\cL
_{M_l}}{{\mathsf y}\in\cL_{M_r}}
\frl^\go_{{\mathsf x}}\mathfrak{t}_{{\mathsf y}-{\mathsf x}}^{\theta
_{\mathsf x}\go}(\frr_{N-M_{r} }^{\theta_{\mathsf y}\go}
- \bfr_{N-M_{r} }).\nonumber
\end{eqnarray}
By the mass separation property \eqref{eq-SepMass}, the annealed
point-to-plane functions~$\bfl_{M}$ and $\bfr_{M}$ have exponentially
decaying tails, and in particular both
are summable.
Define $\mathbf{c}_r \df\sum_M\bfr_{M } <\infty$. The following
theorem is proved in
Section \ref{sub:ProofthmDN}.\vspace*{-2pt}
\begin{thm}
\label{thm:DN}
For every $\lambda>\lambda_0$, there exist $\beta_0 =\beta_0 (\gl
,d)$ and $p_\infty>0$ such that, if Assumption \ref{assA} holds {with $p\leq
p_\infty$}, then, for every $\beta\in[0 ,\beta_0 )$:
\begin{longlist}[(2)]
\item[(1)] The second term on the right-hand side of \eqref{eq:dNDecomp}
converges to zero
$\bbP$-a.s.\ and in $L^2 (\Omega)$.
\item[(2)] The first term on the right-hand side of \eqref{eq:dNDecomp}
converges to
%
%e2.18 ###
\begin{equation}
\label{eq:DNlim}
\frac{\mathbf{c}_r}{\mu} \sum_{{\mathsf x}}\frl_{{\mathsf x}}^\go
\frs
^{\theta_{\mathsf x}\go},
\end{equation}
$\bbP$-a.s.\ and in $L^2 (\Omega)$.\vadjust{\goodbreak}
\end{longlist}
\end{thm}

Consequently, \eqref{eq:limit} of Theorem \ref{Thm:limit} follows with
\[
{\mathfrak{d}^\go}=
\lim_{N\to\infty} \frac{\frD_N^\go}{{\mathbf D}_N^\go} =
\mathbf{c}_r \sum_{{\mathsf x}}\frl_{{\mathsf x}}^\go
\frs^{\theta_{\mathsf x}\go} .
\]
Positivity of $\mathfrak{d}^\go$ on the event $\{0 \in
\operatorname{Cl}_\infty(V)\}$ is established in the
concluding Section~\ref{sub:positivity}.\vspace*{-2pt}

%s3 ###
\section{Proofs}\vspace*{-2pt}
%s3.1 ###
\subsection{The key computation}
\label{ssec-key}
Below, we formulate the key statement, essential for all our results in
this paper. It heavily relies on the assumptions of weak disorder.
We relegate the proof of Proposition \ref{lem:LtwoBound} to the
concluding Section \ref{sec:ProofKey}.\vspace*{-2pt}
\begin{prop}
\label{lem:LtwoBound}
For every $\lambda>\lambda_0$, there exist $\beta_0 =\beta_0 (\gl
,d)$ and \mbox{$p_\infty>0$} such that, if Assumption \ref{assA} holds with $p\leq
p_\infty$, then, for every $\beta\in[0 ,\beta_0 )$,
%
%e3.1 ###
\begin{equation}
\label{eq:key1}
\sup_{N\geq1} \bbE \biggl[ \sum_{{\mathsf x}\in\cH_N^-} \sum
_{{\mathsf y}\in\cH_K^+} \mathfrak{t}^\go_{\mathsf x}(\frq^{\theta
_{\mathsf x}\go}_{{\mathsf y}-{\mathsf x}} - \bfq_{{\mathsf
y}-{\mathsf x}}) g({\mathsf y})
 \biggr]^2 \leqs(K+1)^{1-d/2} \|g\|_{\infty}^2,
\end{equation}
uniformly in $K\geq0$ and in bounded functions $g$ on $\bbZ^{d+1}$.

Furthermore,
%
%e3.2 ###
\begin{equation}
\label{eq:key2}
\bbE \biggl[ \sum_{{\mathsf x}\in\cH_K^-} \sum_{{\mathsf y}\in\cH_K^+}
\mathfrak{t}^\go_{\mathsf x}(\frq^{\theta_{\mathsf x}\go
}_{{\mathsf y}-{\mathsf x}} - \bfq_{{\mathsf y}-{\mathsf x}})
g({\mathsf y})  \biggr]^2
\leqs(K+1)^{ -d/2} \|g\|_{\infty}^2,
\end{equation}
uniformly in $K\geq0$ and in bounded functions $g$ on $\bbZ^{d+1}$.
Similarly,
%
%e3.3 ###
\begin{equation}
\label{eq:key3}
\bbE \biggl[
\sum_{{\mathsf x}, {\mathsf y}\in\cH_K^-}\  \sum_{{\mathsf z}\in
\cH_K^+}
\frl^\go_{{\mathsf x}}\mathfrak{t}_{{\mathsf y}-{\mathsf x}}^{\theta
_{\mathsf x}\go}(\frr_{{\mathsf z}- {\mathsf y}}^{\theta_{\mathsf
y}\go}
- \bfr_{{\mathsf z}-{\mathsf y}})g ({\mathsf z}) \biggr]^2 \leqs
(K+1)^{ -d/2}\|g\|_{\infty}^2
\end{equation}
also uniformly in $K\geq0$ and in bounded functions $g$ on $\bbZ^{d+1}$.\vspace*{-2pt}
\end{prop}

%s3.2 ###
\subsection{\texorpdfstring{Proof of Theorem \protect\ref{thm:L2}(1)}
{Proof of Theorem 2.3(1)}}
\label{sub:ProofthmL21}
Recall that
%
%e3.4 ###
\begin{equation}\label{eq:Difference}
\mathfrak{t}_N^\go- \frac1\mu\frs_N^\go
=
\frac1{2\pi i}\int_{\bbS_\rho} \frac{(1 + \hat\Psi^\go
(u))\Delta(u)}{u^{N+1}}\,\dd u
\end{equation}
for each $\rho\in(0,1 )$.
We are going to show that\vspace*{-2pt}
\begin{lem}
\label{lem:Parseval}
\eqref{eq:Difference} still holds at $\rho=1$ and
$\hat\Psi^\go( e^{i\theta}) \in L^2 (\Omega\times[0,2\pi])$.\vspace*{-2pt}
\end{lem}

In particular,
$\hat\Psi^\go( e^{i\theta}) \in L^2 ( [0,2\pi])$ $\bbP
$-a.s.\
Consequently, the right-hand~si\-de of \eqref{eq:stdifference} is $\bbP
$-a.s.\ equal to the $N$th Fourier coefficient of $(1 + \hat\Psi^\go
( e^{i\theta}))\Delta( e^{i\theta})$. Therefore, by
Parseval's theorem,
\[
\bbE\sum_N  \biggl( \mathfrak{t}_N^\go- \frac1\mu\frs
_N^\go \biggr)^2 = \frac1{2\pi}
\bbE\int_0^{2\pi}  \bigl| \bigl(1 + \hat\Psi^\go( e^{i\theta
})\bigr)\Delta( e^{i\theta})
 \bigr|^2 \,\dd\theta  < \infty.\vadjust{\goodbreak}
\]
It thus follows from Fubini's theorem that
\[
\lim_{N\to\infty} \biggl( \mathfrak{t}_N^\go- \frac
1{\mu} \frs_N^\go \biggr) = 0,
\]
$\bbP$-a.s.\ and in $L^2 (\Omega)$.

It remains to prove Lemma \ref{lem:Parseval}. First of all, $\hat\Psi
^\go( e^{i\theta})$
can be rewritten as
\[
\hat\Psi^\go( e^{i\theta}) = \sum_{{\mathsf x},{\mathsf y}}
\mathfrak{t}_{\mathsf x}^\go(\frq_{{\mathsf y}-{\mathsf x}}^{\theta
_{\mathsf x}
\go} - \bfq_{{\mathsf y}-{\mathsf x}}) e^{i\theta  {\mathsf
y}^{\parallel}}.
\]
Applying Proposition \ref{lem:LtwoBound} with $K=0$, $N=\infty$ and
$g (y) =  e^{i\theta(\vec{e}_{d+1}, {\mathsf y})_{d+1}}$, we
conclude that
\[
\sup_\theta\bbE ( \hat\Psi^\go( e^{i\theta}) )^2
\leqs1 ,
\]
and hence $\Psi^\go(  e^{i\theta} ) \in L_2 (\Omega\times
[0,2\pi])$ indeed.

In a completely similar fashion, one concludes from Proposition \ref
{lem:LtwoBound} that
%
%e3.5 ###
\begin{equation}\label{eq:bigK}
\lim_{K\to\infty} \sup_{|u|\leq1}\bbE \Biggl( \sum_{M=K}^\infty
u^M \psi_M^\go \Biggr)^2 = 0,
\end{equation}
where $\{\psi_N^\go\}$ are the expansion coefficients of $\hat
\Psi^\go(u) = \sum_M u^M \psi_M^\go$, that is, explicitly,
\[
\psi_M^\go= \sum_{{\mathsf y}\in\cL_M}\sum_{{\mathsf
x}}
{\mathfrak t}^\go_{\mathsf x}
(\frq_{{\mathsf y}-{\mathsf x}}^{\theta_{\mathsf x}
\go} - \bfq_{{\mathsf y}-{\mathsf x}}).
\]
Obviously, for each $K$ fixed,
\[
\lim_{\rho\to1} \sum_{N=1}^K  ( \rho e^{i\theta}
)^N \psi_N^\go=
\sum_{N=1}^K  (  e^{i\theta} )^N \psi_N^\go
\]
in $L^2 (\Omega\times[0,2\pi])$. In view of \eqref{eq:bigK},
the latter implies that
\[
\lim_{\rho\to1}\bbE\int_0^{2\pi}  \bigl( \hat\Psi^\go(\rho
 e^{i\theta})
-\hat\Psi^\go(  e^{i\theta})  \bigr)^2\,\dd\theta= 0.
\]
As a result one can indeed pass to the limit $\rho\to1$ in \eqref
{eq:Difference}.\vspace*{-3pt}

%s3.3 ###
\subsection{\texorpdfstring{Proof of Theorem \protect\ref{thm:L2}(2)}
{Proof of Theorem 2.3(2)}}
\label{sub:ProofthmL22}
Let $\cF_N$ be the $\sigma$-algebra generated by $\{V_{\mathsf
x}\}
_{{\mathsf x}\in\cH_{N+1}^-}$, and
let us introduce
\begin{eqnarray*}
\cA^\go_N
&\df&1 + \sum_{M=0}^{N-1} \sum_{{\mathsf x}\in\cL_M}
\sft_{\mathsf x}^\go ( \frq^{\theta_{\mathsf x}\go}_{1,\infty
} -
1 ),\\[-2pt]
\cB^\go_N
&\df&\sum_{M=0}^{N-1} \sum_{{\mathsf x}\in\cL_M}\sft_{\mathsf
x}^\go (
\frq^{\theta_{\mathsf x}\go}_{N-M +1, \infty} - \bfq_{N-M+1
,\infty
} ),\\
\cC^\go_N
&\df&\bbE(\cA^\go_N | \cF_N).
\end{eqnarray*}
We can then express $\frs^\go_N$ as
\[
\frs^\go_N = \cC^\go_N + (\cA^\go_N-\cC^\go_N) - \cB^\go_N.
\]
The $\bbP$-a.s.\ and $L^2 (\Omega)$ convergence in \eqref
{eq:sNslimit} follows from the next two lemmas, since they imply that,
$\bbP$-a.s.\ and in $L^2 (\Omega)$, $\cB^\go_N$ and $\cA^\go
_N-\cC^\go_N$ tend to~$0$, while $\cC^\go_N$ converges to $\frs
^\go$.

\begin{lem}
\label{lem:martingale}
For every $\lambda>\lambda_0$, there exist $\beta_0 >0$ and $\hat
p_\infty>0$ such that, if Assumption \ref{assA} holds with $p \leq\hat
p_\infty$, then, for each $\beta\in[0 ,\beta_0 ]$, the sequence
$\{\cC_N\}$ is an $L^2$-bounded martingale.
\end{lem}

\begin{lem}
\label{lem:compensator}
For every $\lambda>\lambda_0$, there exist $\beta_0>0$ and $\hat
p_\infty>0$ such that, if Assumption \ref{assA} holds with $p \leq\hat
p_\infty$, then
\[
\sum_N\bbE(\cB^\go_N)^2 <\infty \quad \mbox{and}\quad
\sum_N \bbE(\cA^\go_N - \cC^\go_N )^2 <\infty
\]
for each $\beta\in[0 ,\beta_0 ]$.
\end{lem}

\begin{pf*}{Proof of Lemma \ref{lem:martingale}}
The fact that $\cC^\go_N$ is a martingale is straightforward: for any
$N$ and each ${\mathsf x}\in\cL_{N}$,
\[
\bbE( \cC_{N+1}^\go |  \cF_N ) = \bbE( \cA_N^\go |  \cF_N )
+ \sum_{{\mathsf x}\in\cL_N} \bbE \bigl( \mathfrak{t}^\go
_{\mathsf x}(\frq^{\theta_{\mathsf x}\go}_{1, \infty} - 1)|\cF_N
 \bigr) ,
\]
and
\[
\bbE \bigl( \mathfrak{t}^\go_{\mathsf x}(\frq^{\theta_{\mathsf
x}\go}_{1, \infty} - 1)|\cF_N  \bigr) = \mathfrak{t}^\go
_{\mathsf x}\bbE(\frq^{\theta_{\mathsf x}\go}_{1, \infty} - 1) = 0,
\]
since ${\mathsf x}\in\cL_N$.

It remains to check that $\cC^\go_N$ is $L^2(\Omega)$-bounded. We
first deduce from Jensen's inequality that
$\bbE(\cC^\go_N)^2 \leq\bbE(\cA^\go_N)^2$. However,
uniform $L^2 (\Omega)$-boun\-dedness of the latter quantities follows immediately
from Proposition \ref{lem:LtwoBound} by taking $K=0$ and $g\equiv1$.
\end{pf*}

\begin{pf*}{Proof of Lemma \ref{lem:compensator}}
Note that
\[
\cA^\go_N -\cC^\go_N
=
\sum_{{\mathsf x}\in\cH_N^-}\sum_{{\mathsf y}\in\cH
_N^+}
{\mathfrak t}^\go_{\mathsf x} \bigl( \frq_{{\mathsf y}- {\mathsf
x}}^{\theta
_{\mathsf x}\go} -
\bbE ( \frq_{{\mathsf y}- {\mathsf x}}^{\theta_{\mathsf x}\go}
 | \cF_N ) \bigr)
 \]
 and
 \[
\cB^\go_N = \sum_{{\mathsf x}\in\cH_N^-}\sum_{{\mathsf y}\in\cH
_N^+}\mathfrak{t}^\go_{\mathsf x} ( \frq_{{\mathsf y}- {\mathsf
x}}^{\theta_{\mathsf x}\go} - \bfq_{{\mathsf y}- {\mathsf x}}
) .
\]
Thus, $\cA^\go_N -\cC^\go_N$ and $\cB^\go_N$ have very similar
forms. In fact,
\[
\bbE(\cA^\go_N -\cC^\go_N)^2 \leq4 \bbE(\cB_N^\go)
^2 .
\]
On the other hand, taking $g\equiv1$ in the second of the statements of
Proposition \ref{lem:LtwoBound}, we readily conclude that
$\sum_N \bbE(\cB^\go_N)^2 <\infty$.\vadjust{\goodbreak}~%
\end{pf*}

%s3.4 ###
\subsection{\texorpdfstring{Proof of Theorem \protect\ref{thm:DN}}
{Proof of Theorem 2.4}}
\label{sub:ProofthmDN}

The claim (2) of the theorem follows from the \mbox{$\bbP$-a.s.} and
$L^2 (\Omega)$
convergence to $\frs^\go$ in \eqref{eq:sNslimit} and from the fact that
\[
\bbE\frl^\go_{\mathsf x}= \bfl_{\mathsf x}\leq{e }^{-\nu
|{\mathsf x}|}\1
_{\{{\mathsf x}\in\cone\}} .
\]
The first claim (1) follows by an application of \eqref{eq:key3}
with $K=N$ and
$g ({\mathsf z}) = \1_{\{{\mathsf z}\in\cL_N\}}$.

%s3.5 ###
\subsection{\texorpdfstring{Proof of Theorem \protect\ref{Thm:diffusive}}
{Proof of Theorem B}}
\label{sec:ProofThmDiffusive}
Let $f$ be a bounded continuous function on~$\bbR^d$. Using \eqref
{eq:Sinai}, we can write, for any $K\geq0$,
%
%e3.6 ###
\begin{eqnarray}
\label{eq:twosums}
\sum_{{\mathsf z}\in\cL_N}\mathfrak{t}_{\mathsf z}^\go f
\biggl(\frac{{\mathsf z}^\perp}{\sqrt{N}}  \biggr) &=&
\sum_{{\mathsf x}\in\cH_N^-}\sum_{{\mathsf y}}\mathfrak
{t}_{\mathsf x}^\go
(\frq^{\theta_{\mathsf x}\go}_{{\mathsf y}- {\mathsf x}} - \bfq
_{{\mathsf y}- {\mathsf x}})\sum_{{\mathsf z}\in\cL_N}
{\mathbf{ t}} _{{\mathsf z}-{\mathsf y}} f \biggl(\frac{{\mathsf
z}^\perp}{\sqrt{N}} \biggr) \nonumber \\
&=&
\sum_{{\mathsf x}\in\cH_N^-}\sum_{{\mathsf y}\in\cH
_K^-}
{\mathfrak t}_{\mathsf x}^\go
(\frq^{\theta_{\mathsf x}\go}_{{\mathsf y}- {\mathsf x}} - \bfq
_{{\mathsf y}- {\mathsf x}})\sum_{{\mathsf z}\in\cL_N}
{\mathbf{ t}} _{{\mathsf z}-{\mathsf y}} f \biggl(\frac{{\mathsf
z}^\perp}{\sqrt{N}} \biggr) \\
&&{} +
\sum_{{\mathsf x}\in\cH_N^-}\sum_{{\mathsf y}\in\cH
_{K-1}^+}
{\mathfrak t}_{\mathsf x}^\go
(\frq^{\theta_{\mathsf x}\go}_{{\mathsf y}- {\mathsf x}} - \bfq
_{{\mathsf y}- {\mathsf x}})\sum_{{\mathsf z}\in\cL_N}
{\mathbf{ t}} _{{\mathsf z}-{\mathsf y}} f \biggl(\frac{{\mathsf
z}^\perp}{\sqrt{N}} \biggr) .\nonumber
\end{eqnarray}
Choosing $K = K (N) = N^\epsilon$ for some $\epsilon<1/2$ and setting
\[
g({\mathsf y}) = g_N ({\mathsf y})= \sum_{{\mathsf z}\in\cL_N}
{\mathbf{ t}} _{{\mathsf z}-{\mathsf y}} f \biggl(\frac{{\mathsf
z}^\perp}{\sqrt{N}} \biggr),
\]
we can infer from Proposition \ref{lem:LtwoBound} that the second sum
on the right-hand side of \eqref{eq:twosums} converges to zero in
$L^2(\Omega)$. As for the first sum on the right-hand side of \eqref
{eq:twosums}, it follows from the annealed central limit theorem (and
the continuity of $f$) that
\begin{eqnarray*}
&&\lim_{N\to\infty}\max_{{\mathsf y}\in\cH_{N^\epsilon}^-\cap
\cone}
 \biggl|
\sum_{{\mathsf z}\in\cL_N}
{\mathbf{ t}} _{{\mathsf z}-{\mathsf y}} f\biggl(\frac{{\mathsf z}^\perp
}{\sqrt{N}}\biggr)\\
&&\hphantom{\lim_{N\to\infty}\max_{{\mathsf y}\in\cH_{N^\epsilon}^-\cap
\cone}
 \biggl|}{}-
\frac1{\mu}
\frac{1}{\sqrt{\operatorname{det}(2\pi\Sigma)}}\int_{\bbR^d} f({\mathsf x})
 e^{-\gfrac12 (\Sigma^{-1}{\mathsf x},{\mathsf x})}
\,\dd{\mathsf x} \biggr| = 0.
\end{eqnarray*}
By another application of Proposition \ref{lem:LtwoBound}, this time
with $K=0$ and
\begin{eqnarray*}
g ({\mathsf y}) &=&  \biggl(
\sum_{{\mathsf z}\in\cL_N} {\mathbf{ t}} _{{\mathsf z}-{\mathsf y}}
f \biggl(\frac{{\mathsf z}^\perp}{\sqrt{N}} \biggr) -
\frac1{\mu} \frac{1}{\sqrt{\operatorname{det}(2\pi\Sigma)}}\int_{\bbR^d}
f({\mathsf x})
 e^{-\gfrac12 (\Sigma^{-1}{\mathsf x},{\mathsf x})}
\,\dd{\mathsf x} \biggr)\\
&& {}\times\1_{\{{\mathsf y}\in\cH_{N^\epsilon}^-\cap
\cone\}},
\end{eqnarray*}
we conclude, in view of Theorem \ref{thm:L2}, that the first sum
in \eqref{eq:twosums} converges in $L^2(\Omega)$ to
\[
\frac{\frs^\go}{\mu} \frac{1}{\sqrt{\operatorname{det}(2\pi\Sigma
)}}\int_{\bbR^d} f({\mathsf x})
 e^{-\gfrac12 (\Sigma^{-1}{\mathsf x},{\mathsf x})}
\,\dd{\mathsf x}=
\frac{\mathfrak{t}^\go}{\sqrt{\operatorname{det}(2\pi\Sigma
)}}\int_{\bbR^d} f({\mathsf x})
 e^{-\gfrac12 (\Sigma^{-1}{\mathsf x},{\mathsf x})}
\,\dd{\mathsf x}.
\]

\subsubsection*{Extension to the full $\cD_N$-ensemble}
Proceeding as in Section \ref{sub:DN}, we conclude that, for any
bounded continuous $f$, the series
\[
\sum_{{\mathsf z}\in\cL_N}\mathfrak{d}_{\mathsf z}^\go f
\biggl(\frac{{\mathsf z}^\perp}{\sqrt{N}}  \biggr)
\]
converges in $L^2(\Omega)$ to
\[
\frac{
\mathbf{c}_r \sum_{{\mathsf x}}\frl_{{\mathsf x}}^\go\frs^{\theta
_{\mathsf x}\go}
}{\mu\sqrt{\operatorname{det}(2\pi\Sigma)}}\int_{\bbR^d} f({\mathsf x})
 e^{-\gfrac12 (\Sigma^{-1}{\mathsf x},{\mathsf x})}
\,\dd{\mathsf x}.
\]
Together with \eqref{eq:DNlim}, this implies \eqref{eq:flim}.

\subsubsection*{Local limit description}
As in \cite{Sinai}, equations \eqref{eq:Sinai} and \eqref{eq:Sinaib}
suggest the following
quenched Ornstein--Zernike asymptotics for $\mathfrak{t}_{\mathsf
x}^\go$ (as inherited from the annealed OZ-asymptotics of
${\mathbf{ t}} _{\mathsf x}$ in \cite{IV-annealed}): Given
${\mathsf x}\in\bbZ^{d+1}$, let $\hat{\theta}_{\mathsf x}\go$ be the
reflection with respect to the hyperplane $\cL_N$ of the shifted
environment $\theta_{\mathsf x}\go$. In other words, $\hat{\theta
}_{\mathsf x}\go$ is the environment as seen backwards from ${\mathsf
x}$. Of course,
the reflected environment has the very same averaged polymer
connectivity functions. We conjecture that
%
%e3.7 ###
\begin{equation}\label{eq:OZRandom}
\frac{\mathfrak{t}_{\mathsf x}^\go}{{\mathbf{ t}} _{\mathsf x}} = (
1+\frs^\go)( 1+\frs^{\hat\theta_{\mathsf x}\go})\bigl(1+\so\bigr).
\end{equation}
Clearly, the strength of the above conjecture depends on what is meant
by~$\so$ in \eqref{eq:OZRandom}. A $\bbP$-a.s.\ statement would be a
refinement of a $\bbP$-a.s.\ CLT, which is, as we already mentioned,
an open problem by itself. Weaker statements, on the other hand, are
feasible via an appropriate refinement of Proposition~\ref{lem:LtwoBound}.

%s4 ###
\section{\texorpdfstring{$L^2 (\Omega)$ estimates at weak disorder}
{L2 (Omega) estimates at weak disorder}}
\label{sec:ProofKey}
%s4.1 ###
\subsection{Preliminaries}
\label{sub:Prilim}
Our proof of Proposition \ref{lem:LtwoBound} is based on a comparison
with weakly interacting random walks on $\bbZ^d$. The bottom line is
that, under Assumption \ref{assA}, transience wins over attraction. From a
technical point of view the approach is similar to \cite{CIL}.

Since, in all the estimates below, only the supremum norm of $g$ in
Proposition~\ref{lem:LtwoBound} would matter, we can assume, without
loss of generality, that $g\equiv1$.
It is convenient to use the alternative notation
\[
\frq_{{\mathsf x},{\mathsf u}}^\go
\df\sum_{\gamma\in\cT_{{\mathsf x},{\mathsf u}}^0} W^\go_{\lambda
,\beta
}(\gamma)
= \sum_{\gamma\in\cT_{{\mathsf u}-{\mathsf x}}^0} W^{\theta
_{\mathsf x}\go
}_{\lambda,\beta}(\gamma) ,
\]
and $\bfq_{{\mathsf x},{\mathsf u}} = \bfq_{{\mathsf u}-{\mathsf x}}
\df\bbE\frq
_{{\mathsf x},{\mathsf u}}^\go$.
Above, $\cT_{{\mathsf u}}^0$ is the set of irreducible cone-confined paths
from $0$ to ${\mathsf u}$ and $\cT_{{\mathsf x},{\mathsf u}}^0 \df
{\mathsf x}+ \cT
_{{\mathsf u}-{\mathsf x}}^0 $.

Given ${\mathsf x}$ and ${\mathsf u}$, we define the diamond shape
\[
D({\mathsf x},{\mathsf u}) \df({\mathsf x}+\cone)\cap(
{\mathsf u}-\cone).\vadjust{\goodbreak}
\]
By construction, any path $\gamma\in\cT_{{\mathsf x},{\mathsf u}}$ satisfies
$\gamma\subset D({\mathsf x},{\mathsf u})$. Hence, $\frq^\go
_{{\mathsf x},{\mathsf u}}$ only depends on the environment inside
$D({\mathsf x},{\mathsf u})$.

%f2 ###
\begin{figure}

\includegraphics{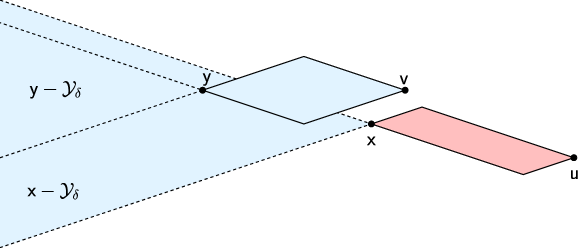}%
\vspace*{-3pt}
\caption{With ${\mathsf x},{\mathsf y},{\mathsf u},{\mathsf v}$ as in
the picture, the
paths contributing to $\mathfrak{t}_{\mathsf x}^\go
\mathfrak{t}_{\mathsf y}^\go( \frq^\go_{{\mathsf y},{\mathsf v}} -
\bfq_{{\mathsf y},{\mathsf v}})$ lie inside the blue region, while the paths
contributing to $\frq^\go_{{\mathsf x},{\mathsf u}} - \bfq_{{\mathsf
x},{\mathsf u}}$ lie inside the red region. These two quantities are thus
independent (w.r.t.\ the disorder).}
\label{fig:indep}
\vspace*{-3pt}
\end{figure}
Here is a useful observation (see Figure \ref{fig:indep}): If $D (
{\mathsf x},{\mathsf u})\cap D ({\mathsf y},{\mathsf v})=
\emptyset$, then
\[
\bbE \{ \mathfrak{t}_{\mathsf x}^\go(\frq^\go
_{{\mathsf x},{\mathsf u}} - \bfq_{{\mathsf x},{\mathsf u}})
\mathfrak{t}_{\mathsf y}^\go(\frq^\go_{{\mathsf y},{\mathsf v}}
- \bfq_{{\mathsf y},{\mathsf v}}) \} = 0.
\]
Indeed, unless ${\mathsf x}= {\mathsf y}$, it is always true that
either $D (
{\mathsf x},{\mathsf u})\cap({\mathsf y}-\cone)= \emptyset$ or
$D ({\mathsf y},{\mathsf v})\cap({\mathsf x}-\cone)=
\emptyset$. If,
in addition, the diamond shapes do not intersect, then in the former
case $(\frq^\go_{{\mathsf x},{\mathsf u}} - \bfq_{{\mathsf
x},{\mathsf u}})$
is independent of $\mathfrak{t}_{\mathsf x}^\go
{\mathfrak t}_{\mathsf y}^\go(\frq^\go_{{\mathsf y},{\mathsf v}}
- \bfq
_{{\mathsf y},{\mathsf v}})$, and similarly for the latter case.

Consequently, neglecting nonpositive terms, we obtain
\begin{eqnarray*}
&&\bbE \biggl\{
\sum_{{\mathsf x}\in\cH_N^-} \sum_{{\mathsf y}\in\cH_K^+}
{\mathfrak t}^\go_{\mathsf x}(\frq^{\theta_{\mathsf x}\go
}_{{\mathsf y}-{\mathsf x}} - \bfq_{{\mathsf y}-{\mathsf x}})  \biggr\}
^2\\[-3pt]
&&\qquad \leq\mathop{\sum_{{\mathsf x},{\mathsf y}\in\cH_N^-}}_{{\mathsf
u},{\mathsf v}\in
\cH_K^+}
\bbE\{
\mathfrak{t}_{\mathsf x}^\go\frq^\go_{{\mathsf x},{\mathsf u}}
\mathfrak{t}_{\mathsf y}^\go\frq^\go_{{\mathsf y},{\mathsf v}} +
\mathfrak{t}_{\mathsf x}^\go\bfq_{{\mathsf x},{\mathsf u}}
\mathfrak{t}_{\mathsf y}^\go\bfq_{{\mathsf y},{\mathsf v}} \}
\1_{\{ D ({\mathsf x},{\mathsf u})\cap D ({\mathsf y},{\mathsf
v})\neq
\emptyset\}} .
\end{eqnarray*}
Now, it follows from the attractiveness \eqref{eq-attract} of the
interaction that
\[
\bbE\{
\mathfrak{t}_{\mathsf x}^\go\frq^\go_{{\mathsf x},{\mathsf u}}
\mathfrak{t}_{\mathsf y}^\go\frq^\go_{{\mathsf y},{\mathsf v}} \}
\geq
\bbE\{\mathfrak{t}_{\mathsf x}^\go\bfq_{{\mathsf x},{\mathsf u}}
\mathfrak{t}_{\mathsf y}^\go\bfq_{{\mathsf y},{\mathsf v}} \},
\]
and thus
\begin{eqnarray}
\label{eq:BoundDiamond}
&&\bbE \biggl\{
\sum_{{\mathsf x}\in\cH_N^-} \sum_{{\mathsf y}\in\cH_K^+}
{\mathfrak t}^\go_{\mathsf x}(\frq^{\theta_{\mathsf x}\go
}_{{\mathsf y}-{\mathsf x}} - \bfq_{{\mathsf y}-{\mathsf x}})
\biggr\}^2\nonumber\\[-10pt]\\[-10pt]
&&\qquad \leq2  \mathop{\sum_{{\mathsf x},{\mathsf y}\in\cH_N^-}}_{
{\mathsf u},{\mathsf v}\in\cH_K^+}
\bbE\{
\mathfrak{t}_{\mathsf x}^\go\frq^\go_{{\mathsf x},{\mathsf u}}
\mathfrak{t}_{\mathsf y}^\go\frq^\go_{{\mathsf y},{\mathsf v}} \}
\1_{\{ D ({\mathsf x},{\mathsf u})\cap D ({\mathsf y},{\mathsf
v})\neq
\emptyset\}} .\nonumber
\end{eqnarray}
The latter expression sets up the stage for an analysis in terms of weakly
interacting random walks.\vadjust{\goodbreak}

%s4.2 ###
\subsection{Weakly interacting random walks}
\label{sub:weak}
Let $\Prw$ be the path measure of a random walk on $\bbZ^{d+1}$ whose
independent steps are distributed
according to $\{\bfq_\ell\}$. We shall
use notation $\underline{X} = (X_0 , X_1, \dots)$ for the path
of this random walk. Let us
say that $({\mathsf x},{\mathsf u})\in\underline{\sfX}$ if there
exists $n$
such that
$\sfX_n = {\mathsf x}$ and $\sfX_{n+1} ={\mathsf u}$.
In this way,
$\Prw(({\mathsf x},{\mathsf u})\in\underline{\sfX})
={
{\mathbf t}} _{\mathsf x}\bfq_{{\mathsf x},{\mathsf u}}$.
Let also $\Prwt$ be the product measure for a couple of such random walks.

Given a path $\underline{{\mathsf x}} = (0={\mathsf x}_0 ,{\mathsf
x}_1 ,{\mathsf x}_2, \dots)$, we define the random functionals
\[
\cQ_n^\go(\underline{{\mathsf x}}) \df
\prod_{i=1}^n
\frq_{{\mathsf x}_{i-1} ,{\mathsf x}_i}^\go.
\]
Note that $\bbE\cQ_n^\go(\underline{{\mathsf x}}) = \Prw(
\underline
{\sfX}_n =
\underline{{\mathsf x}})$, where the event $\{\underline{\sfX}_n=
\underline{{\mathsf x}}\}$ means that the first $n$ steps of
$\underline{\sfX}$ are given by the corresponding steps of
$\underline{{\mathsf x}}$.

Consider now two admissible trajectories $\underline{{\mathsf x}}$ and
$\underline{{\mathsf y}}$. For any $n\in\bbN$, we define the diamond sausage
$D(\underline{{\mathsf x}}_n )$ around the first $n$ steps of
$\underline
{{\mathsf x}}$ by
\[
D(\underline{{\mathsf x}}_n ) \df\bigcup_1^n D ({\mathsf x}_{i-1}
,{\mathsf x}_i
) .
\]
By definition, $D(\underline{{\mathsf x}} ) \df D(\underline{{\mathsf
x}}_\infty) $. If
$ D(\underline{{\mathsf x}}_n )\cap D(\underline{{\mathsf y}}_m) =
\emptyset
$, then
\[
\bbE\cQ_n^\go(\underline{{\mathsf x}})\cQ_m^\go(\underline
{{\mathsf y}}) =
\Prwt(\underline{\sfX}_n =\underline{{\mathsf x}} ; \underline
{\sfY}_m =
\underline{{\mathsf y}} ).
\]
If, however, the above diamond sausages intersect, then, by the positive
association \eqref{eq-attract} of one-dimensional random variables,
%
%e4.1 ###
\begin{equation}
\label{eq:PosAs}
\bbE\cQ_n^\go(\underline{{\mathsf x}})\cQ_m^\go(\underline
{{\mathsf y}})
\geq
\Prwt(\underline{\sfX}_n =\underline{{\mathsf x}} ; \underline
{\sfY}_m =
\underline{{\mathsf y}}),
\end{equation}
which means that the random weights $\cQ^\go$ produce attraction
between the
two paths. In particular, all terms which contribute to the right-hand side
of \eqref{eq:BoundDiamond} satisfy
\[
\bbE\{
\mathfrak{t}_{\mathsf x}^\go\frq^\go_{{\mathsf x},{\mathsf u}}
\mathfrak{t}_{\mathsf y}^\go\frq^\go_{{\mathsf y},{\mathsf v}} \}
\geq\Prwt\bigl(({\mathsf x},{\mathsf u})\in\underline{\sfX} ;
({\mathsf y},{\mathsf v})\in\underline{\sfY}
\bigr).
\]
Let now $\underline{{\mathsf x}}$ and $\underline{{\mathsf y}}$ be two
infinite admissible
paths. We define the corresponding diamond intersection number
\[
\#( \underline{{\mathsf x}} ,\underline{{\mathsf y}})\df
\# \{ (k,\ell) \dvtx  D({\mathsf x}_{k-1} ,{\mathsf x}_k)\cap
D({\mathsf y}_{\ell
-1} ,{\mathsf y}_\ell)\neq
\emptyset \} .
\]
Let also $\cE$ be the event that there exist $k,\ell$ such that
$D({\mathsf x}_{k-1} ,{\mathsf x}_k)\cap D({\mathsf y}_{\ell-1}
,\allowbreak{\mathsf y}_\ell
)\neq\emptyset$, ${\mathsf x}_{k-1},{\mathsf y}_{\ell-1}\in\cH
_N^-$ and
${\mathsf x}_k,{\mathsf y}_\ell\in\cH_K^+$.
Expanding $\mathfrak{t}_{\mathsf x}^\go$ and $
{\mathfrak t}_{\mathsf y}^\go$ as in the first line of
\eqref{eq:Sinai1}, we infer that
the sum on the right-hand side of \eqref{eq:BoundDiamond} is bounded
above by
\begin{eqnarray}
\label{eq:XYDsum}
&&\sum_{\underline{{\mathsf x}} ,\underline{{\mathsf y}}}
\bbE\cQ^\go(\underline{{\mathsf x}}) \cQ^\go(\underline{{\mathsf y}})
 \#( \underline{{\mathsf x}} ,\underline{{\mathsf y}} ) \1_{\cE
}(\underline
{{\mathsf x}} ,\underline{{\mathsf y}})\nonumber\\[-8pt]\\[-8pt]
&&\qquad \df\lim_{m\to\infty}\lim_{n\to\infty}
\sum_{\underline{{\mathsf x}} ,\underline{{\mathsf y}}}
\bbE\cQ^\go_{m} (\underline{{\mathsf x}}) \cQ^\go_{n} (\underline
{{\mathsf y}})
 \#( \underline{{\mathsf x}} ,\underline{{\mathsf y}} ) \1_{\cE
}(\underline
{{\mathsf x}} ,\underline{{\mathsf y}}) .\nonumber
\end{eqnarray}
Existence of the above limit follows by monotonicity from \eqref{eq:PosAs}.
We thus obtain
%
%e4.2 ###
\begin{equation}
\label{eq:RWBound}\qquad
\bbE \biggl\{
\sum_{{\mathsf x}\in\cH_N^-} \sum_{{\mathsf y}\in\cH_K^+}
{\mathfrak t}^\go_{\mathsf x}(\frq^{\theta_{\mathsf x}\go
}_{{\mathsf y}-{\mathsf x}} - \bfq_{{\mathsf y}-{\mathsf x}})  \biggr\}^2
\leq
2\sum_{\underline{{\mathsf x}} ,\underline{{\mathsf y}}}
\bbE \{
\cQ^\go(\underline{{\mathsf x}})\cQ^\go(\underline{{\mathsf y}})
\# (\underline{{\mathsf x}}, \underline{{\mathsf y}} ) \1_{\cE}
 \}.
\end{equation}
Of course, in order to apply
the latter upper bound
one needs to control the statistics of
% may be huge or even infinite if, e.g.,
$\# (\underline{{\mathsf x}}, \underline{{\mathsf y}} )$.
% = \infty$.
The point is that,
under Assumption \ref{assA}, the $\cQ^\go$-induced interaction between the
paths $\underline{\sfX} $ and $\underline{\sfY}$ is so weak that it
does not
destroy transient behavior.
This phenomenon is stated in Lemma \ref{lem:tworatio} below, in a form
which happens to be particularly convenient for the latter use.

Given $t$, ${\mathsf u}_0 ,{\mathsf v}_0\in\cL_0$ and ${\mathsf u}_1
,{\mathsf v}_1\in\cL_t$ consider two pieces
$\underline{{\mathsf x}}_n$ and $\underline{{\mathsf y}}_m$ of admissible trajectories (assuming that they exist): $\underline
{{\mathsf x}} = (
{\mathsf u}_0 = {\mathsf x}_0, \dots,{\mathsf x}_n ={\mathsf u}_1,
\dots)$ from
${\mathsf u}_0$ to ${\mathsf u}_1$, and
$\underline{{\mathsf y}} = ({\mathsf v}_0 = {\mathsf y}_0, \dots
,{\mathsf y}_m
={\mathsf v}_1 ,\dots)$
from ${\mathsf v}_0$ to ${\mathsf v}_1$.
\begin{lem}
\label{lem:tworatio}
Once $\lambda>\lambda_0$ is fixed, for
every $\eta>0$ there exists $\beta_0 >0$ and $p_{\infty}>0$ such that
%
%e4.3 ###
\begin{equation}
\label{eq:tworatio}\qquad
\bbE\cQ^\go_n(\underline{{\mathsf x}} )\cQ^\go_m (\underline
{{\mathsf y}} )
\leq
\exp\bigl\{\tfrac12\eta t\1_{\{
D (\underline{{\mathsf x}}_n )\cap D(\underline{{\mathsf y}}_m )\neq
\emptyset
\}}\bigr\}\Prwt(\underline{\sfX}_n =
\underline{{\mathsf x}}, \underline{\sfY}_m = \underline{{\mathsf
y}}),
\end{equation}
uniformly in $\beta\in[0,\beta_0)$, provided that Assumption \ref{assA} is
satisfied with $p<p_{\infty}$.
The inequality \eqref{eq:tworatio} holds simultaneously for all
$t$, ${\mathsf u}_0 ,{\mathsf v}_0\in\cL_0$ and ${\mathsf u}_1
,{\mathsf v}_1\in\cL
_t$ and the corresponding admissible
trajectories $\underline{{\mathsf x}} $, $\underline{{\mathsf y}}$.
\end{lem}

\begin{pf}%{Proof of Lemma \ref{lem:tworatio}}
The left-hand side of \eqref{eq:tworatio} equals
to $\Prwt(\underline{\sfX}_n =
\underline{{\mathsf x}}, \underline{\sfY}_m = \underline{{\mathsf
y}})$
%$1$
whenever
$D (\underline{{\mathsf x}}_n )\cap D(\underline{{\mathsf y}}_m ) =
\emptyset$.
Indeed, in such a situation, $\cQ^\go_n(\underline{{\mathsf x}})$
and $\cQ^\go_m (\underline{{\mathsf y}})$ are independent.

We proceed to consider the case when
$D (\underline{{\mathsf x}}_n )\cap D(\underline{{\mathsf y}}_m )\neq
\emptyset$.
Let us say that a path $\gamma\in\cT_{{\mathsf u}_0 ,{\mathsf u}_1}$ is
compatible with
$\underline{{\mathsf x}}_n;\ \gamma\sim\underline{{\mathsf x}}_n$, if
$\underline{{\mathsf x}}_n\setminus\{{\mathsf x}_0 ,{\mathsf
x}_n\}$ is
precisely the collection of all
the cone points of $\gamma$. Similarly for $\gamma^\prime\sim
\underline{{\mathsf y}}_m$.
The left-hand side in \eqref{eq:tworatio} is
\[
 e^{2t\xi}\sumtwo{\gamma\sim\underline{{\mathsf x}}_n}
{\gamma'\sim\underline{{\mathsf y}}_n}\bbE
W_{\gl,\gb}^\go(\gamma) W_{\gl,\gb}^\go(\gamma' ) =
\sumtwo{\gamma\sim\underline{{\mathsf x}}_n}{\gamma'\sim
\underline
{{\mathsf y}}_m}
\exp\{2t\xi-\lambda(\vert\gamma\vert +|\gamma' |)-\Phi
_\gb(\gamma, \gamma' )
\},
\]
where the annealed interaction potential
$\Phi_\gb(\gamma, \gamma' )$
is given by
\[
\Phi_\gb(\gamma, \gamma' ) =
\sum_{\sfw\in\bbZ^{d+1}} \phi_\gb(\ell_{\gamma\cup\gamma'}
(\sfw))
\qquad \mbox{with }
\phi_\gb(\ell) \df-\log\bbE e^{-\ell V^\go} .
\]
Above, $\ell_{\gamma\cup\gamma'} (\sfw)$ is the total combined local
time of the couple $( \gamma,\gamma' )$ in $\sfw$. Therefore,
ignoring the interaction, one derives the following upper bound:
\[
\bbE\cQ^\go_n(\underline{{\mathsf x}} )\cQ^\go_m (\underline
{{\mathsf y}} )\leq
\sumtwo{\gamma\sim\underline{{\mathsf x}}_n}{\gamma'\sim
\underline
{{\mathsf y}}_n}
\exp\{2t\xi-\lambda(\vert\gamma\vert +|\gamma' |)\},
\]
that is, in terms of the corresponding expression for the simple
symmetric random walk on
$\bbZ^{d+1}$ with the constant killing rate $\gl-\gl_0 = \gl-\log
(2d ) > 0$.

Similarly,
\[
\Prwt(\underline{\sfX}_n =
\underline{{\mathsf x}}, \underline{\sfY}_m = \underline{{\mathsf
y}})=
\sumtwo{\gamma\sim\underline{{\mathsf x}}_n}{\gamma'\sim
\underline
{{\mathsf y}}_m}
\exp\{2t\xi-\lambda(\vert\gamma\vert +|\gamma' | )-\Phi
_\gb(\gamma) -
\Phi_\gb(\gamma' ) \}.
\]
The function $\phi_\gb$ is subadditive \cite{Flury,IV-annealed}.
Consequently $\phi_\gb(\ell)
\leq\ell\phi_\gb(1)$.
We conclude that the following lower bound on
$\Prwt(\underline{\sfX}_n =
\underline{{\mathsf x}}, \underline{\sfY}_m = \underline{{\mathsf
y}})$
holds for any $c >0$:
\[
 e^{-c t\phi_\gb(1) }
\sumtwo{\gamma\sim\underline{{\mathsf x}}_n}{\gamma'\sim
\underline
{{\mathsf y}}_m}
\exp\{2t\xi-\lambda(\vert\gamma\vert +|\gamma' |)\}
1_{\{\vert\gamma\vert +|\gamma' | \leq ct\}} .
\]
Recall that $t$ is the horizontal span of both $\gamma$ and $\gamma'$
and that
$\gl>\gl_0 = \log(2d)$ is fixed. Thus, as directly follows from the
properties
of the simple random walk on $\bbZ^{d+1}$ subject to a constant
killing potential $ \gl-\gl_0$, there
exists $\epsilon= \epsilon(c)$, tending to zero as $c\to\infty$,
such that
\[
\sumtwo{\gamma\sim\underline{{\mathsf x}}_n}{\gamma'\sim
\underline
{{\mathsf y}}_m}
\exp\{-\lambda(\vert\gamma\vert +|\gamma' | )\}1_{\{
\vert\gamma\vert +|\gamma' | \leq ct\}} \geq
\bigl(1- \epsilon(c) \bigr)
\sumtwo{\gamma\sim\underline{{\mathsf x}}_n}{\gamma'\sim
\underline
{{\mathsf y}}_m}
\exp\{-\lambda(\vert\gamma\vert +|\gamma' | )\}.
\]
Altogether, we conclude that, for any $c >0$,
\[
\frac{\bbE\cQ^\go_n(\underline{{\mathsf x}} )\cQ^\go_m
(\underline
{{\mathsf y}} )}
{ \Prwt(\underline{\sfX}_n =
\underline{{\mathsf x}}, \underline{\sfY}_m = \underline{{\mathsf
y}})}
\leq
\exp\bigl\{ct\phi_\gb(1) - \log\bigl(1 - \epsilon(c )\bigr)\bigr\}.
\]
In its turn, the smallness of $\phi_\gb(1)$ is controlled through
\[
\lim_{\beta\to0} \phi_\gb(1) = -\log(1- p) .
\]
Consequently, the claim of the Lemma follows first by taking $c$
sufficiently large and then
by choosing $\beta$ and $p$ appropriately small.
\end{pf}

%s4.3 ###
\subsection{Upper bounds in terms of synchronized random walks}
\label{sub:synchro}
\begin{figure}

\includegraphics{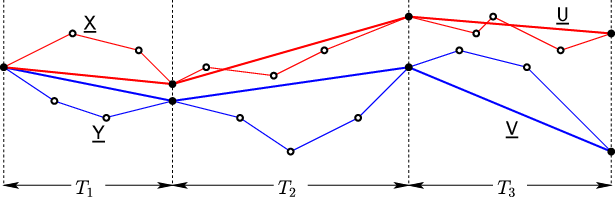}

\caption{The $\underline{\sfX}$, $\underline{\sfY}$ and $(
\underline{\sfU}, \underline{\sfV})$ random walks.}
\label{fig:XYUV}
\end{figure}

Let us explain how Lemma \ref{lem:tworatio} is put to work in order to
control \eqref{eq:XYDsum}.
At this stage, it happens to be convenient to synchronize the two
trajectories $\underline{\sfX}$ and $\underline{\sfY}$, by
expressing all the above quantities in terms of another induced $\bbZ
\times\bbZ^d\times\bbZ^d$-valued random walk $(\underline{\sfU
}, \underline{\sfV} )$: Let $\underline{{\mathsf x}}$ and
$\underline
{{\mathsf y}}$ be realizations of $\underline{\sfX}$ and $\underline
{\sfY}$. Let us label all the $\cL_n$-hyperplanes which are
simultaneously hit by both the $\underline{{\mathsf x}}$ and~$\underline
{{\mathsf y}}$ trajectories as $n_1, n_2, \ldots,$ with ${\mathsf
u}_1, {\mathsf u}_2, \dots$ and ${\mathsf v}_1 , {\mathsf v}_2 ,\dots
$ the corresponding
hitting points (see Figure \ref{fig:XYUV}).
Then the induced trajectory of $(\underline{\sfU}, \underline
{\sfV})$ is
$(\underline{{\mathsf u}}, \underline{{\mathsf v}})$.
We denote by $t_1, t_2, \ldots$ the horizontal spans of the steps of
$(\underline{{\mathsf u}}, \underline{{\mathsf v}})$.
We shall use\vadjust{\goodbreak} $\PRWS$ for the path measure of $(\underline{\sfU},
\underline{\sfV})$. The distribution of a single step under~$\PRWS
$ is given by
\begin{eqnarray*}
\PRWS({\mathsf u}, {\mathsf v})& =& \PRWS(t ,u, v) \\
&=&
\mathop{\sum_{n=1}}_{m=1}^t  \mathop{\sum_{0<t_1<\dots<t_n =t}}_{
0 <s_1 <\dots< s_m = t}
\mathop{\sum_{{\mathsf x}_i\in\cL_{t_i}}}_{{\mathsf y}_j\in\cL_{s_j}}
\prod_1^n\bfq_{{\mathsf x}_i - {\mathsf x}_{i-1}}\prod_1^m\bfq
_{{\mathsf y}_i -
{\mathsf y}_{i-1}}
\prodtwo{0 <i<n}{0<j <m}\1_{\{ t_i\neq s_j\}} ,
\end{eqnarray*}
where we have set ${\mathsf x}_0 = {\mathsf y}_0 = 0$ and ${\mathsf
x}_n = {\mathsf u},  {\mathsf y}_m ={\mathsf v}$.
Alternatively,
\[
\PRWS({\mathsf u}, {\mathsf v}) = \PRWS(t ,u, v) = \Prwt \bigl(
T(\underline
{\sfX} ,\underline{\sfY})=t ;
{\mathsf u}\in\operatorname{Range}(\underline{\sfX} ); {\mathsf v}\in\operatorname{Range}(\underline{\sfY} ) \bigr) ,
\]
where
\[
T(\underline{\sfX} ,\underline{\sfY}) \df
\inf \{ n\dvtx  \operatorname{Range}(\underline{\sfX})\cap\cL_n
\neq\emptyset  \mbox{ and }
\operatorname{Range}(\underline{\sfY})\cap\cL_n\neq\emptyset \}
\]
is the (random) horizontal span of a step of the $(\underline{\sfU
}, \underline{\sfV})$-random walk.
In view of the uniform exponential tails of $\{\bfq_N\}$, there
exists $\kappa=\kappa(\lambda) >0$ such that
%
%e4.4 ###
\begin{equation}\label{eq:TkDecat}
\PRWS(T >\ell)\leqs e^{-\kappa\ell} ,
\end{equation}
uniformly in $l$ and in $\beta\geq0$.

%f3 ###

Let us go back to \eqref{eq:XYDsum}.
The i.i.d.\ horizontal spans of $(\underline{\sfU}, \underline
{\sfV})$-steps will be denoted by
$T_1, T_2, \dots.$
To ease notation, set $D_k (\underline{{\mathsf u}}) \df D({\mathsf
u}_k ,{\mathsf u}_{k+1})$ and similarly for $D_k (\underline{{\mathsf
v}})$. Obviously,
if a pair
$(\underline{{\mathsf x}}, \underline{{\mathsf y}})$ of $( \underline
{\sfX
}, \underline{\sfY})$-paths is compatible with a synchronized
$(\underline{\sfU}, \underline{\sfV})$-path
$(\underline{{\mathsf u}}, \underline{{\mathsf v}})$;
$(\underline{{\mathsf x}}, \underline{{\mathsf y}})\sim(\underline
{{\mathsf u}}, \underline{{\mathsf v}})$,
then
\[
\# (\underline{{\mathsf x}}, \underline{{\mathsf y}}) \leq\sum_{k}
T_k
\1_{\{ D_k (\underline{{\mathsf u}})\cap D_k (\underline{{\mathsf
v}}) \neq
\emptyset\}}.
\]
By Lemma \ref{lem:tworatio}, once $\lambda>\lambda_0$ is fixed,
for every $\eta>0$ there exist $\beta_0 >0$ and $p_{\infty}>0$ such that\footnote{Strictly speaking, the inequality makes sense for restrictions
to any finite number of steps of the
$(\underline{{\mathsf u}}, \underline{{\mathsf v}})$-trajectory.}
\[
\sum_{(\underline{{\mathsf x}}, \underline{{\mathsf y}})\sim
(\underline
{{\mathsf u}}, \underline{{\mathsf v}})}
\bbE
\cQ^\go(\underline{{\mathsf x}})\cQ^\go(\underline{{\mathsf y}})
\leq\exp \biggl\{ \frac12 \eta\sum_k T_k
\1_{\{ D_k (\underline{{\mathsf u}})\cap D_k (\underline{{\mathsf v}})
\neq\emptyset\}} \biggr\}
\PRWS(\underline{{\mathsf u}}, \underline{{\mathsf v}}).
\]
Therefore, since $x e^x \leq e^{2x}$ for all $x\geq0$,
\eqref{eq:RWBound} implies that
%
%e4.5 ###
\begin{eqnarray}\label{eq:RWBoundUV}
&&\bbE \biggl\{
\sum_{{\mathsf x}\in\cH_N^-} \sum_{{\mathsf y}\in\cH_K^+}
{\mathfrak t}^\go_{\mathsf x}(\frq^{\theta_{\mathsf x}\go
}_{{\mathsf y}-{\mathsf x}} - \bfq_{{\mathsf y}-{\mathsf x}})
\biggr\}^2\nonumber\\[-8pt]\\[-8pt]
&&\qquad\leq
\frac{2}{\eta}
\ERWS
\exp \biggl\{ \eta\sum_k T_k
\1_{\{ D_k (\underline{\sfU})\cap D_k (\underline{\sfV}) \neq
\emptyset\}}  \biggr\} \1_{\widehat\cE} ,\nonumber
\end{eqnarray}
where $\widehat\cE$ is the analog of $\cE$ for the synchronized
random walks, that is,
\[
\widehat\cE\df \{\exists k   \dvtx
D_k (\underline{\sfU})\cap D_k (\underline{\sfV}) \neq\emptyset;
\sfU_{k},\sfV_k \in\cH_N^- \mbox{ and } \sfU_{k+1} ,\sfV
_{k+1}\in\cH_K^+ \} .
\]
Of course $\cE\subset\widehat\cE$, in the sense that if $\cE$
holds for
$( \underline{{\mathsf x}}, \underline{{\mathsf y}} )$, then
$\widehat\cE$ also
holds for the synchronized $( \underline{{\mathsf u}}, \underline
{{\mathsf v}} )$ path.

Let us now bound the expectation in the right-hand side
of \eqref{eq:RWBoundUV}, uniformly in $\eta$ sufficiently small.
Let $\sfZ_k\df\sfU_k-\sfV_k$, and notice that there exists a
constant $\alpha=\alpha(d,\delta)$ such that
\[
\exp \bigl\{ \eta T_k
\1_{\{ D_k (\underline{\sfU})\cap D_k (\underline{\sfV}) \neq
\emptyset\}}
 \bigr\}
\leq
\exp \bigl\{ \eta T_k
\1_{\{ T_k > \alpha\|\sfZ_{k-1}^{\perp}\|\}}
 \bigr\}.
\]
Writing $\exp \{\eta T_k\1_{\{ T_k > \alpha\|\sfZ
_{k-1}^{\perp}\|\}} \} =  (( e^{\eta
T_k} - 1 )\1_{\{ T_k > \alpha\|\sfZ_{k-1}^{ \perp
}\|\}} + 1  )$ and expanding, we obtain
\[
\exp \Biggl\{ \sum_{k=1}^M \eta T_k
\1_{\{ T_k > \alpha\|\sfZ_{k-1}^{\perp}\|\}}
 \Biggr\}
=
\sum_{A\subset\{1,\ldots,M\}}
\prod_{k\in A}  (  e^{\eta T_k}- 1  ) \1_{\{ T_k >
\alpha\|\sfZ_{k-1}^{\perp}\|\}}
.
\]
Since $( e^{\eta T_k}-1)/( e^\eta-1) \leq T_k  e^{\eta
T_k}$, we can bound the right-hand side from above by
%
%e4.6 ###
\begin{equation}
\label{eq:intermediate}
\sum_{n\geq0} ( e^\eta-1)^n \mathop{\sum_{A\subset\{
1,\ldots,M\}}}_{|A|=n} \prod_{k\in A}
T_k e^{\eta T_k} \1_{\{ T_k > \alpha\|\sfZ
_{k-1}^{\perp}\|\}}
\1_{\widehat{\cE}} .
\end{equation}
Let us write $A=\{a_1,a_2,\ldots,a_n\}$, with $a_1 < a_2 <\cdots
<a_n$, and let us set $a_0=0$. We are going to split the trajectories
into $n$ ``bubbles,'' the $i$th bubble being composed of the steps
$\sfZ_{a_{i-1}+1},\ldots,\sfZ_{a_i}$. The horizontal span $B_i$ of
the $i$th bubble is thus
\[
B_i \df\sum_{k=a_{i-1}+1}^{a_i} T_k,\qquad  1\leq i \leq n.
\]

%s4.4 ###
\subsection{\texorpdfstring{Proof of Proposition \protect\ref{lem:LtwoBound}}
{Proof of Proposition 3.1}}
\label{sub:Key}
We only prove \eqref{eq:key1} and \eqref{eq:key2}, the third claim,
\eqref{eq:key3}, being a variant of the latter.

We first prove \eqref{eq:key1}. In this case, we only retain from the
event $\widehat\cE$ the constraint that $\sum_i B_i > K$.
More precisely, we bound above the\vadjust{\goodbreak} $\ERWS$-expectation of the sum in
\eqref{eq:intermediate} by
%
%e4.7 ###
\begin{equation}
\label{eq:intermediate1}\qquad
\ERWS\sum_{n\geq1} ( e^\eta-1)^{n-1}
\sum_{
|
A|=n
%}
}
\1_{\{ \sum_i B_i > K\}}
\prod_{k\in A}
T_k e^{\eta T_k} \1_{\{ T_k > \alpha\|\sfZ
_{k-1}^{\perp}\|\}} .
\end{equation}
Therefore, by the Markov property, \eqref{eq:RWBoundUV} implies
\begin{eqnarray} \label{eq:bubbles}
&&\sup_N\bbE \biggl\{
\sum_{{\mathsf x}\in\cH_N^-} \sum_{{\mathsf y}\in\cH_K^+}
{\mathfrak t}^\go_{\mathsf x}(\frq^{\theta_{\mathsf x}\go
}_{{\mathsf y}-{\mathsf x}} - \bfq_{{\mathsf y}-{\mathsf x}})
\biggr\}^2\nonumber\\[-8pt]\\[-8pt]
&&\qquad \leq
\frac2\eta\sum_{n\geq1} ( e^\eta-1)^{n-1} \mathop{\sum_{
B_1,\ldots,B_n}}_{\sum_i B_i > K} \prod_{i=1}^n I(B_i),\nonumber
\end{eqnarray}
where, for $B\in\bbN$,
\[
I(B) \df\sup_{{\mathsf z}\in\bbZ^d} \sum_{m=1}^B \ERWS \Biggl(
{
T_m  e^{\eta T_m }   ;
}
\sum_{k=1}^{m}
T_k = B, T_m > \alpha\|Z_{m-1}^{\perp}\|
| \sfZ_0^{\perp} = {\mathsf z} \Biggr).
\]
We need a reasonable upper bound on the latter quantities.
Recall that we can choose $\eta$ as small as we wish.
Observe first that \eqref{eq:TkDecat} and a standard large deviation
estimate imply the existence of $\epsilon>0$ and $c>0$ such that,
uniformly in $B\in\bbN$,
\[
\sup_{{\mathsf z}\in\bbZ^d} \sum_{m=1}^{\epsilon B} \ERWS \Biggl(
{
T_m  e^{\eta T_m }
}
  ;
\sum_{k=1}^{m} T_k = B \Big| \sfZ_0^{\perp} =
{\mathsf z} \Biggr) \leqs e^{-cB}.
\]
On the other hand, relying again on \eqref{eq:TkDecat} and using the
local limit theorem for i.i.d.\ random variables with exponential
tails, we obtain that
\begin{eqnarray}\label{eq:Bcomp}
&&\sup_{{\mathsf z}\in\bbZ^d} \sum_{m=\epsilon B}^B \ERWS \Biggl(
{
T_m  e^{\eta T_m }
};
\sum_{k=1}^{m} T_k = B, T_m > \alpha\|Z_{m-1}^{
\perp}\| | \sfZ
_0^{\perp} = {\mathsf z} \Biggr)\nonumber \\
&&\qquad \leqs
\sum_{t\geq1} t e^{-(\nu-\eta) t} \sup_{{\mathsf z}\in\bbZ^d}
\sum_{m=\epsilon B}^B \PRWS ( \|Z_{m-1}^{
\perp}\| < t/\alpha| \sfZ
_0^{\perp} = {\mathsf z} )\nonumber\\[-8pt]\\[-8pt]
&&\qquad \leqs
\sum_{t\geq1}  e^{-(\nu-\eta) t} \frac{t^{d+1}}{B^{d/2}}\nonumber \\
&&\qquad \leqs
B^{-d/2}.\nonumber
\end{eqnarray}
We therefore conclude that, for any $B\in\bbN$,
\[
I(B) \leqs B^{-\gfrac d2}.
\]
Let us now use this bound to control the right-hand side of \eqref{eq:bubbles}.
For~fi\-xed~$n$, let $L=\sum_{i=1}^n B_i$; then there must be an index
$j$ such that $\prod_{i=1}^n I(B_i) \leqs(n/L)^{\gfrac d2} \prod
_{i\neq j} I(B_i)$. Therefore, choosing $\eta$ small enough, we have
%
%e4.8 ###
\begin{eqnarray}\label{eq:Lcomp}
\nonumber
&&\sum_{n\geq1} ( e^\eta-1)^{n-1} \mathop{\sum_{B_1,\ldots
,B_n}}_{\sum_i B_i > K} \prod_{i=1}^n I(B_i)\\[-2pt]
&&\qquad \leqs
\sum_{L>K} L^{-\gfrac d2} \sum_{n\geq0} n^{1+\gfrac d2}  \Biggl(
( e^\eta-1)\sum_{B\geq1} I(B)  \Biggr)^n \\[-2pt]
&&\qquad \leqs
\sum_{L>K} L^{-\gfrac d2}
\leqs
(1+K)^{1-\gfrac d2}.\nonumber
\end{eqnarray}

Let us now turn to the proof of \eqref{eq:key2}.

We proceed to bound the right-hand side of
\eqref{eq:key2} in terms
of the synchronized random walks $\underline\sfU$ and
$\underline\sfV$.
As before, $\sfZ_k = \sfU_k -\sfV_k$. Let
$j_0$ be such that $\sfZ_{j_0-1}\in\cH_K^-$ and $\sfZ_{j_0}\in\cH
_K^+$. We need to
derive a bound on
\[
\ERWS\exp \biggl\{ \sum_k \eta T_k
\1_{\{ T_k > \alpha\|\sfZ_{k-1}^{\perp}\|\}}
 \biggr\}
\1_{\{ D_{j_0-1}(
\underline\sfU)\cap D_{j_0 -1}(\underline\sfV)\neq\emptyset\}} .
\]

Expanding as in \eqref{eq:intermediate1}, we may restrict attention to
sets $A$ which
contain an element $a_{i_0}$ such that $a_{i_0}=j_0$. This implies
that, if $\sum_{i=1}^{i_0} B_i = K+t$, the excess $t$ must be entirely
due to the $j_0$th step of $\sfZ$. In particular, this quantity
has exponential tails, and, following the derivation of
\eqref{eq:Bcomp},
\[
I (B_{i_0}) \leqs e^{-(\nu-\eta)t}B_{i_0}^{-\gfrac{d}{2}} .
\]
We can thus write, proceeding as in \eqref{eq:Lcomp},
\begin{eqnarray*}
&&\bbE \biggl\{
\sum_{{\mathsf x}\in\cH_K^-} \sum_{{\mathsf y}\in\cH_K^+}
{\mathfrak t}^\go_{\mathsf x}(\frq^{\theta_{\mathsf x}\go
}_{{\mathsf y}-{\mathsf x}} - \bfq_{{\mathsf y}-{\mathsf x}})
\biggr\}^2\\[-2pt]
&&\qquad \leq
\frac2\eta\sum_{t\geq1}\sum_{i_0\geq1} \sum_{n\geq0} ( e^\eta-1)^{
{
n+i_0 -1
}
} \mathop{\sum_{B_1,\ldots,B_{n+i_0}}}_{\sum_{i=1}^{i_0} B_i =
K+t} \prod_{i=1}^{n+i_0} I(B_i)\\[-2pt]
&&\qquad \leqs
\sum_{t\geq1} e^{-(\nu- \eta) t} \sum_{i_0\geq1} ( e^\eta
-1)^{i_0} \mathop{\sum_{B_1,\ldots,B_{i_0}}}_{\sum_{i=1}^{i_0} B_i
= K+t} \prod_{i=1}^{i_0} B_i^{-\gfrac d2}\\[-2pt]
&&\qquad \leqs
\sum_{t\geq1} e^{-(\nu-\eta) t} (K+t)^{-\gfrac d2}\\[-2pt]
&&\qquad \leqs
(1+K)^{-\gfrac d2}.
\end{eqnarray*}

\begin{rem}
\label{rem:Trans}
The above computations readily imply the following: Let ${\mathsf
u},{\mathsf v}\in\cL_0$ and let
$\PRWS_{{\mathsf u},{\mathsf v}}$ be the distribution of the
synchronized $(\underline{\sfU}, \underline{\sfV})$ random
walk starting
from $({\mathsf u},{\mathsf v})$. Then, under Assumption \ref{assA},
%
%e4.9 ###
\begin{equation}
\label{eq:Trans}
\ERWS_{{\mathsf u},{\mathsf v}}\exp \biggl\{ \eta\sum_k T_k\1
_{\{ D_k (
\underline{\sfU})
\cap D_k ( \underline{\sfV})\neq\emptyset\}} \biggr\} \leqs1 ,
\end{equation}
uniformly in ${\mathsf u},{\mathsf v}$ and in all $\eta$ sufficiently small.
\end{rem}

%s4.5 ###
\subsection{\texorpdfstring{Positivity of $\mathfrak{d}^\go$ on the event $\{0\in \operatorname{Cl}_\infty(V)\}$}
{Positivity of d omega on the event \{0 in Cl infinity(V)\}}}
\label{sub:positivity}

Let $0\in\operatorname{Cl}_\infty(V)$. Then $\mathfrak{d}^\go>0$
if there exists ${\mathsf x}= (x ,t)$ such that
$\mathfrak{d}^{\theta_{\mathsf x}\go} >0$. Indeed, such
${\mathsf x}$ should necessarily
satisfy ${\mathsf x}\in\operatorname{Cl}_\infty(V)$. Hence, there exists a
nearest-neighbor
\textit{finite} path $\gamma= (\gamma(0), \dots, \gamma(n))$
from $0$ to
${\mathsf x}$ such that $\gamma(l)\in\operatorname{Cl}_\infty(V)$ for all
$l=0, \dots, n$ and, consequently, such that
$W_{\lambda,\gb}^\go(\gamma) >0$. However,
\[
\frD_N^\go\geq W_{\lambda,\gb}^\go(\gamma) \frD_{N- t}^{\theta
_{\mathsf x}\go}.
\]
It follows that
\[
\liminf_{N\to\infty}
 e^{N\xi} \frD_N^\go\geqs e^{t \xi}
W_{\lambda,\gb}^\go(\gamma) \frs^{\theta_{\mathsf x}\go}.
\]
It remains to show that
\[
\bbP ( \exists  {\mathsf x}\dvtx  \frs^{\theta_{\mathsf x}\go} >0 )= 1.
\]
In fact, an ostensibly stronger claim holds:

\begin{lem}
Under conditions of Theorem \ref{thm:L2},
\label{lem:exists}
\[
\bbP ( \exists  {\mathsf x}\in\cL_0 :  \frs^{\theta
_{\mathsf x}\go
} >0 ) = 1.
\]
\end{lem}

\begin{pf}%{Proof of Lemma \ref{lem:exists}}
The proof is by the second moment method, and based on
$L^2$-esti\-mates at weak disorder as developed in the
preceding subsection. Let $B_n\subset\cL_0$ be the $d$-dimensional
lattice box of side-length $n$,
\[
B_n \df \bigl\{ {\mathsf x}= (x_1, \dots,x_d ,0)  \dvtx  x_l\in\{0,
\dots, n-1\}
\mbox{ for } l=1,\dots, d \bigr\} .
\]
By Theorem \ref{thm:L2}, $\bbE\frs^{\theta_{\mathsf x}\go} \equiv1$.
We claim that the variance
%
%e4.10 ###
\begin{equation}
\label{eq:BoxVar}
\operatorname{\mathbb{V}ar}\biggl(\frac1{n^d}\sum_{{\mathsf x}\in B_n}  ( \frs^{\theta
_{\mathsf x}\go} - 1 )\biggr)
\leqs\frac1{n^{d/2 -1}} .
\end{equation}
The conclusion of the lemma would then follow by
Chebyshev's estimate and a Borel--Cantelli argument. Now, the estimates
developed in
%Lemma \ref{lem:tworatio} of
the preceding subsections
% \ref{sub:Key}
imply that, under Assumption \ref{assA}, the extra attraction
stemming from integration of the factors $\cQ^\go$ over intersecting diamonds
does not alter the statistical properties of the effective
$d$-dimensional random walks $(\underline{\sfX},
\underline{\sfY})$, or, equivalently,
of the synchronized random walks $(\underline{\sfU}, \underline{\sfV
})$. In particular, for any
${\mathsf x}, {\mathsf y}\in B_n$,
%
%e4.11 ###
\begin{equation}
\label{eq:xyBoxBound}
\qquad  |
\bbE ( \frs^{\theta_{\mathsf x}\go} - 1 )
 ( \frs^{\theta_{\mathsf y}\go} - 1 )  | \leqs
\PRWS_{{\mathsf x},{\mathsf y}} \bigl(D(\underline{\sfU})\cap
D(\underline{\sfV})\neq\emptyset\bigr)
\leqs\frac1{|{\mathsf x}- {\mathsf y}|^{d/2 - 1}} .
\end{equation}
Indeed, the second inequality above is straightforward.
As for the first inequality in \eqref{eq:xyBoxBound}, proceeding as in
the proof of \eqref{eq:RWBoundUV}, we infer that
\[
 |
\bbE ( \frs^{\theta_{\mathsf x}\go} - 1 )
 ( \frs^{\theta_{\mathsf y}\go} - 1 )  | \leqs
\ERWS
\exp \biggl\{ \eta\sum_k T_k
\1_{\{ D_k (\underline{\sfU})\cap D_k (\underline{\sfV}) \neq
\emptyset\}}  \biggr\} \1_{\{ D(\underline{\sfU})\cap D(\underline
{\sfV})\neq\emptyset\}} .
\]
By the strong Markov property and in view of \eqref{eq:Trans},
\[
\ERWS
\exp \biggl\{ \eta\sum_k T_k
\1_{\{ D_k (\underline{\sfU})\cap D_k (\underline{\sfV}) \neq
\emptyset\}}  \biggr\} \1_{\{ D(\underline{\sfU})\cap D(\underline
{\sfV})\neq\emptyset\}}
\leqs
\PRWS_{{\mathsf x},{\mathsf y}} \bigl(D(\underline{\sfU})\cap
D(\underline
{\sfV})\neq\emptyset\bigr),
\]
and \eqref{eq:xyBoxBound} follows.

The variance decay estimate
\eqref{eq:BoxVar} is a direct consequence of \eqref{eq:xyBoxBound}
\[
\operatorname{\mathbb{V}ar}\biggl(\frac1{n^d}\sum_{{\mathsf x}\in B_n}  ( \frs^{\theta
_{\mathsf x}\go} - 1 )\biggr)
\leqs\frac1{n^{2d}}\cdot n^d \cdot\sum_{k=1}^n \frac
{k^{d-1}}{k^{d/2 -1}} .
\]
\upqed
\end{pf}

\section*{Acknowledgment}

We thank Francis Comets for pointing out the
similarity between our method and
the expansion used in \cite{Sinai} to treat the directed case.

% imsref loaded by dianan, 2011-02-04 13:36:19
%

\printaddresses

\end{document}